\documentclass[12pt]{amsart}
\usepackage{amssymb,verbatim}
\newtheorem{theorem}{Theorem}[section]
\newtheorem{proposition}[theorem]{Proposition}
\newtheorem{conjecture}[theorem]{Conjecture}
\newtheorem{corollary}[theorem]{Corollary}
\newtheorem{lemma}[theorem]{Lemma}

\theoremstyle{definition}

\newtheorem{remark}[theorem]{Remark}
\newtheorem{example}[theorem]{Example}
\newtheorem{definition}[theorem]{Definition}
\newtheorem{problem}[theorem]{Problem}

\newcounter{romctr}
\renewcommand{\theromctr}{\alph{romctr}}
\newenvironment{romitem}[1]
{\refstepcounter{romctr}%
\begin{itemize}
\item[(\theromctr)] 
#1}
{\end{itemize}}

\newenvironment{smallbmatrix}[1]
{\left[\begin{smallmatrix}
#1}
{\end{smallmatrix}\right]}

\renewcommand{\eqref}[1]{{\rm (\ref{#1})}}

\newcommand{\Hom}{ \operatorname{Hom}}
\newcommand{\LR}{ \operatorname{LR}}
\newcommand{\Horn}{ \operatorname{{\bf H}}}
\newcommand{\lr}[2]{c_{\,#1}^{\,#2}}
\newcommand{\eqbydef}{\stackrel{\rm def}{=\!\!=}}

\begin{document}

{\ }\vspace{-.3in}

\title[Eigenvalues, singular values, and
Littlewood-Richardson Coefficients]
{Eigenvalues, singular values, and\\
Littlewood-Richardson Coefficients}

\date{August 9, 2003}

\author{Sergey Fomin}
\address{Department of Mathematics \\
University of Michigan \\
Ann Arbor, MI 04109
}
\email{fomin@umich.edu}

\author{William Fulton}
\address{Department of Mathematics \\
University of Michigan \\
Ann Arbor, MI 04109
}
\email{wfulton@umich.edu}

\author{Chi-Kwong Li}
\address{Department of Mathematics \\
College of William and Mary \\
Williamsburg, VA 23187
}
\email{ckli@math.wm.edu}

\author{Yiu-Tung Poon}
\address{Department of Mathematics \\
Iowa State University \\
Ames, IA 50011}
\curraddr{Department of Mathematics \\
Hong Kong Institute of Education \\
10 Lo Ping Road, 
Tai Po, N.~T., 
Hong Kong
}
\email{ytpoon@iastate.edu, ytpoon@ied.edu.hk}

\thanks{Partially supported by
NSF grants DMS-0070685, 
DMS-9970435, 
and 
DMS-0071994. 
}

\subjclass{Primary
15A42. 
Secondary
05E15, 
14M15, 
15A18.
}

\begin{abstract}
We characterize the relationship between the singular values of 
a Hermitian (resp., real symmetric, complex symmetric) matrix
and  the singular values of its off-diagonal block.
We also characterize the eigenvalues of a Hermitian (or real
symmetric) matrix $C=A+B$ in terms of the combined list of eigenvalues
of $A$ and~$B$. 
The answers are given by Horn-type linear inequalities. 
The proofs depend on a new inequality among Littlewood-Richardson coefficients.
\end{abstract}

\maketitle

\section*{Introduction}

Let $X$ be the upper right $p$ by $n\!-\!p$ submatrix of an $n$ by $n$
matrix~$Z$, with $2p \leq n$.  
If $Z$ is (complex) Hermitian or real symmetric, the main result of~\cite{LP} 
characterizes the possible eigenvalues of~$Z$ in terms of the
singular values of~$X$.  
In this paper, we provide the analogous characterization 
for the singular values of~$Z$, 
when $Z$ is Hermitian, or complex symmetric, or real symmetric.  
Surprisingly, the possibilities in all three cases are the same.

The answers are given by linear inequalities of Horn type.  
Let $\gamma_1 \geq \cdots \geq \gamma_n$ and
$s_1 \geq \cdots \geq s_p$ be the singular values of $Z$ and~$X$,
respectively. 
We prove that the possible pairs of sequences $(s_k)$ and
$(\gamma_\ell)$ are characterized by the inequalities 
\[
2 \sum_{k \in K} s_k \leq \sum_{i \in I}
\gamma_{2i-1} + \sum_{j\in J} \gamma_{2j}\,,
\]
for all triples $(I,J,K)\in\bigcup_{r\leq p}\LR_r^p$, 
where $\LR_r^p$ denotes the list of triples 
defined inductively by Horn~\cite{H}.
It appears~in relation to a number of problems surveyed
in~\cite{Fu1} (where $\LR_r^p$ is denoted by~$T_r^p$),
including the original Horn's problem of characterizing the 
eigenvalues of a sum $A+B$ of two Hermitian (or real symmetric)
matrices in terms of the eigenvalues of $A$ and~$B$.
The definition of $\LR_r^p$ used in this paper 
describes the triples $(I,J,K)\in \LR_r^p$ as those for which
a certain \emph{Littlewood-Richardson coefficient} does not vanish. 
See Definition~\ref{def:lr-set}.

In Section~\ref{sec:horn-combined} of this paper, 
we solve the following modification of Horn's
problem: what are the possible eigenvalues of $A+B$ given the
combined list of eigenvalues of $A$ and~$B$ (without specifying
which eigenvalues belong to~$A$, and which ones to~$B$)? 

Our proofs depend on Klyachko's celebrated solution of Horn's problem,
on its refinement obtained by Knutson and Tao, 
on the results on eigenvalues and singular
values from \cite{Fu1, Fu2, Kl2, LP}, and 
on a version of the Littlewood-Richardson rule given in~\cite{CL}. 

\section{Singular value inequalities}

\begin{definition}
\label{def:lr-set}
For positive integers $r\leq p$, the set $\LR_r^p$ 
consists of all triples $(I,J,K)$ of subsets of
$\{1, \dots, p\}$ of the same cardinality $r$,
such that the Littlewood-Richardson coefficient 
$\lr{\lambda(I) \, \lambda(J)}{\lambda(K)}$ is positive.
(For general background on Little\-wood-Richard\-son coefficients,
see for example \cite{Fu0} or~\cite{SF-EC2}.) 
Here and in what follows, we use the correspondence 
\begin{equation}
\label{eq:set-to-partition}
I = \{i_1 < i_2 < \cdots < i_r \}\ \ \mapsto\ \ {}
\lambda(I) = (i_r - r, \dots, i_2 -2, i_1 - 1) 
\end{equation}
between $r$-element subsets $I$ 
(always written with this notation in increasing order) 
and integer partitions~$\lambda(I)$ with at most $r$ parts. 
\end{definition}

For any complex $\ell$ by $m$ matrix $X$, there are unitary matrices
$U$ and $V$ 
($\ell$ by $\ell$ and $m$ by $m$ respectively) such that $UXV$ has
nonnegative entries $s_1 \geq s_2 \geq \ldots \geq s_p$
in positions $(1,1), (2,2),\dots,(p,p)$,  where $p = \min(\ell, m)$, 
and all other entries zero.  The real numbers $s_1,
\ldots, s_p$, always written in decreasing order, with multiplicities,
are the \emph{singular values} of~$X$.
The positive singular values are the positive square roots of the
eigenvalues of the positive semidefinite matrix~$X^* X$, where $X^*$
denotes the conjugate transpose of~$X$. 
If $X$ is Hermitian, then its singular values are simply 
the absolute values of its eigenvalues. 

\subsection{Main result}

The following is our main result on singular values. 

\begin{theorem} 
\label{th:main}
Let $p$ and $n$ be positive integers, with $2p \leq n$, and let
$\gamma_1 \geq \cdots \geq \gamma_n$ and
$s_1 \geq \cdots \geq s_p$ be sequences of nonnegative
real numbers.  The following are equivalent:

\begin{romitem}
\label{any-complex}
There exists a complex matrix $Z$ of the form 
$
Z=\begin{smallbmatrix}
* & X \\
Y & * 
\end{smallbmatrix}
$
such that
\begin{itemize}
\item[$\bullet$]
$X$ is $p$ by $n-p$, with singular values $s_1,\dots,s_p$; 
\item[$\bullet$]
$Y$ is $n-p$ by $p$, with singular values $s_1,\dots,s_p$; 
\item[$\bullet$]
$Z$ is $n$ by $n$, with singular values $\gamma_1,\dots,\gamma_n$.  
\end{itemize}
\end{romitem}
\begin{romitem}
\label{real-symmetric}
There exists a real symmetric 
matrix $Z=\begin{smallbmatrix}
* & X \\
X^* & * 
\end{smallbmatrix}$ 
such that
\begin{itemize}
\item[$\bullet$]
$X$ is $p$ by $n-p$, with singular values $s_1,\dots,s_p$; 
\item[$\bullet$]
$Z$ is $n$ by $n$, with singular values $\gamma_1,\dots,\gamma_n$
(thus, eigenvalues $\pm\gamma_1,\dots,\pm\gamma_n$).  
\end{itemize}
%
\end{romitem}
\begin{romitem}
\label{alternating-eigenvalues}
There exists a real symmetric 
matrix $Z=\begin{smallbmatrix}
* & X \\
X^* & * 
\end{smallbmatrix}$ 
such that
\begin{itemize}
\item[$\bullet$]
$X$ is $p$ by $n-p$, with singular values $s_1,\dots,s_p$; 
\item[$\bullet$]
$Z$ is $n$ by $n$, with eigenvalues 
$\gamma_1,-\gamma_2,\gamma_3,-\gamma_4,\dots,(-1)^{n-1}\gamma_n$.
\end{itemize}
\end{romitem}
\begin{romitem}
\label{lr-inequalities}
The numbers $\gamma_1,\dots,\gamma_{2p}$ and $s_1,\dots,s_p$ satisfy
the linear inequalities
\begin{equation}
\label{eq:lr-inequality}
2 \sum_{k \in K} s_k \leq \sum_{i \in I}
\gamma_{2i-1} + \sum_{j\in J} \gamma_{2j},
\end{equation}
for all $r\leq p$ and all triples $(I,J,K)\in\LR_r^p$. 
\end{romitem}
\end{theorem}

The equivalence of \eqref{alternating-eigenvalues} and
\eqref{lr-inequalities} is a direct corollary of the main results in
\cite{Fu2} and~\cite{LP} 
(see Section~\ref{sec:compare-to-LP} below). 
The other 
equivalences are new. 
Conditions \eqref{any-complex}--\eqref{lr-inequalities} are also
equivalent to condition~\eqref{three-hermitian} 
in Corollary~\ref{cor:three-hermitian},
and to condition~\eqref{three-hermitian-combined} in
Corollary~\ref{cor:sv-ev}.

\begin{remark}
In \eqref{any-complex}, 
``complex''
can be replaced by ``complex symmetric''
or ``Hermitian.'' 
Indeed, the resulting statements imply~\eqref{any-complex}
and are implied by \eqref{real-symmetric}.
In these versions, the somewhat unnatural requirement that $X$ and $Y$ have the same
singular values becomes redundant. 
See also Remark~\ref{rem:pool-together}. 
\end{remark}

\begin{remark}
Condition \eqref{lr-inequalities} does not involve the $\gamma_i$
with~$i>2p$. 
Hence in \eqref{any-complex}--\eqref{alternating-eigenvalues},
one can replace 
\emph{``$Z$ is $n$ by $n$, with singular values
  $\gamma_1,\dots,\gamma_n$''} 
by \emph{``$Z$ is $n$ by $n$, with largest $2p$ singular values 
$\gamma_1,\dots,\gamma_{2p}\,$.''} 
\end{remark}

\begin{remark}
In \eqref{lr-inequalities}, it suffices to check the inequalities for
the triples $(I,J,K)$ with $\lr{\lambda(I)\,\lambda(J)}{\lambda(K)}=1$.  
This follows from a theorem of P.~Belkale (see \cite[Proposition 9]{Fu1}).
\end{remark}

\begin{remark}
\label{rem:specify-in-advance}
The matrices $X$ and $Y$ in \eqref{any-complex} can be
specified in advance,
as can the matrix $X$ in \eqref{real-symmetric}--\eqref{alternating-eigenvalues}.  
Indeed, any two matrices with the same singular
values can be trans\-formed into each 
other by multiplying on the left and right by unitary (orthogonal in
the real case) matrices.
On the other hand, for unitary matrices $U_1$, $U_2$, $V_1$, and
$V_2$, the matrix 
\begin{equation} 
\label{eq:block-unitary}
\begin{bmatrix}U_1 & 0 \\ 0 & U_2\end{bmatrix} \cdot \begin{bmatrix}P
      & X \\ Y & Q\end{bmatrix}
\cdot \begin{bmatrix}V_1 & 0 \\ 0 & V_2\end{bmatrix} \, = \,
\begin{bmatrix}U_1PV_1 & U_1XV_2 \\ U_2YV_1 & U_2QV_2\end{bmatrix}
\end{equation}
has the same singular values as $Z=\begin{smallbmatrix}P
      & X \\ Y & Q\end{smallbmatrix}$. 
\end{remark}

\subsection{Comparison with previous results. Examples}
\label{sec:compare-to-LP}

We next summarize the main result of~\cite{LP}, which
includes the main result of~\cite{Fu2}, in 
a form suitable for our purposes.  

\begin{theorem}
\label{th:LiPoon}
Let $2p \leq n$, let $s_1 \geq \ldots \geq s_p \geq 0$, 
and $\lambda_1\geq \ldots \geq \lambda_n$. 
The following are equivalent:
\begin{enumerate}
\item[(i)]
There exists an $n$ by $n$ Hermitian matrix of the form
      $\begin{smallbmatrix}* & X \\ X^* & *\end{smallbmatrix}$
      with eigenvalues
$\lambda_1, \ldots, \lambda_n$, such that
$X$ is $p$ by $n-p$ with singular values $s_1, \ldots, s_p$.
\item[(ii)]
For all $r \leq p$ and $(I,J,K)$ in $\LR_r^p$,
\[
2\sum_{k \in K} s_k \, \leq \,
\sum_{i \in I} \lambda_i - \sum_{j \in J} \lambda_{n+1-j}.
\]
\item[(iii)]
There exist Hermitian $p$ by $p$ matrices $A$, $B$, and $C$ with
eigenvalues $\lambda_1, \ldots, \lambda_p$, $\lambda_{n+1-p}, \ldots,
\lambda_n$, and $s_1, \ldots, s_p$, 
respectively, such that $2C \leq A - B$, i.e., $A - B - 2C$ is
positive semidefinite.
\item[(iv)]  
There exist Hermitian $p$ by $p$ matrices $A$, $B$, and $C$
such that $2C \leq A - B$,
the eigenvalues of $C$ are $s_1, \ldots, s_p$,
and the 
eigenvalues of~$A$ (respectively,~$B$)
listed  in descending order (with multiplicities) form a
subsequence of $\lambda_1, \ldots, \lambda_n$. 
\end{enumerate}
The Hermitian matrices in \textnormal{(i)}, \textnormal{(iii)}, and 
\textnormal{(iv)} can
be taken to be real symmetric matrices, and the matrices $X$ in
\textnormal{(i)} and $C$ in \textnormal{(iii)} and 
\textnormal{(iv)} 
can be specified in advance.
\end{theorem}

The equivalences (i)--(iii) in Theorem~\ref{th:LiPoon} are essentially 
\cite[Theorem~2.3]{LP}. 
The proof of Theorem~\ref{th:LiPoon} is given in Section~\ref{sec:proof-of-th:LiPoon}. 

\pagebreak[2]

\begin{corollary}
\label{cor:three-hermitian}
In Theorem~\ref{th:main}, 
conditions~\eqref{alternating-eigenvalues} and \eqref{lr-inequalities}
are equivalent to each other, and also to 
\begin{romitem}
\label{three-hermitian}
There exist Hermitian $p$ by $p$ matrices $A$, $B$, and~$C$ 
with eigenvalues $\gamma_1, \gamma_3, \ldots, \gamma_{2p-1}$ (of $A$),
$\gamma_2, \gamma_4, \ldots, \gamma_{2p}$ (of $B$), and $s_1, s_2, \ldots,
s_p$ (of $C$), such that $2C \leq A + B$.
\end{romitem}
\end{corollary}

\begin{proof}
Apply the equivalence (i)$\Leftrightarrow$(iii) of Theorem~\ref{th:LiPoon}, 
with 
$(\lambda_1,\lambda_2,\dots,\lambda_{n-1},\lambda_n)
=(\gamma_1,\gamma_3, \dots, -\gamma_4,-\gamma_2)$. 
(Here in the right-hand side, 
the numbers $(-1)^{i-1}\gamma_i$ are listed in decreasing order.) 
\end{proof}

\begin{remark}
Condition \eqref{three-hermitian} can be further relaxed---see
condition~\eqref{three-hermitian-combined} in Corollary~\ref{cor:sv-ev}. 
\end{remark}

\begin{remark}
\label{rem:cones}
Condition \eqref{real-symmetric} of Theorem~\ref{th:main} 
concerns real symmetric (or Hermitian) matrices with
eigenvalues $\pm \gamma_1, \dots, \pm \gamma_n$.
There are $2^n$ choices of signs, and for each choice of signs,
Theorem~\ref{th:LiPoon} describes the relationship between the
possible values of $\gamma_1, \dots, \gamma_n$ 
and $s_1,\dots,s_p$ by a set of linear inequalities~(ii). 
The assertion
$\eqref{real-symmetric}\Leftrightarrow\eqref{alternating-eigenvalues}$
in Theorem~\ref{th:main} can be rephrased as saying that 
the set of inequalities corresponding to the alternating choice of
signs $(\gamma_1,-\gamma_2,\gamma_3,-\gamma_4,\dots)$
is uniformly the weakest.
That is, the inequalities (ii) corresponding to any other choice of signs
imply the inequalities \eqref{eq:lr-inequality}. 
\end{remark}

\begin{remark}
Removing the restrictions concerning matrix $Y$ in part
\eqref{any-complex} of  Theorem~\ref{th:main} 
results in a different (strictly weaker
than~\eqref{lr-inequalities})---and much simpler---set of
inequalities. 
Specifically, it is not hard to prove (cf.\ a more general result of
Thompson's \cite{T}) that, for $2p\leq n$, there exists a
$p$ by $n-p$ matrix with singular values
$s_1 \ge \cdots \ge s_p$ that is a submatrix of an
$n$ by $n$ matrix with singular values
$\gamma_1 \ge \cdots \ge \gamma_n$ if
and only if $\gamma_i \ge s_i$ for $i = 1, \dots, p$.
\end{remark}

\begin{remark}
O'Shea and Sjamaar~\cite{OS} give another polyhedral description of
the singular values appearing in
part~\eqref{real-symmetric} of Theorem~\ref{th:main}. 
However, \cite{OS}~contains no inequalities like~\eqref{eq:lr-inequality}.
We do not know whether the equivalence
\eqref{real-symmetric}$\Leftrightarrow$\eqref{lr-inequalities}
can be deduced~from~\cite{OS}. 
\end{remark}


\begin{example}
Let $p=1$.
The only triple $(I,J,K)\in\LR_1^1$ 
is $(\{1\},\{1\},\{1\})$. \linebreak[2]
The corresponding inequality \eqref{eq:lr-inequality} is 
\begin{equation}
\label{eq:p=1}
2s_1\leq \gamma_1+\gamma_2\,.
\end{equation} 
In the special case $p=1$, $n=2$, the equivalence
\eqref{any-complex}$\Leftrightarrow$\eqref{lr-inequalities} 
in Theorem~\ref{th:main} can be restated as follows: 
the singular values 
$\gamma_1\geq\gamma_2\geq 0$
of a $2$ by $2$ complex matrix of the form 
$\begin{smallbmatrix}
* & x \\ 
y & *
\end{smallbmatrix}$,
where 
$|x|=|y|=s_1\geq 0$ is fixed, satisfy the inequality
\eqref{eq:p=1} 
and no other constraints. 
Equivalently (see Remark~\ref{rem:specify-in-advance}), 
the inequality \eqref{eq:p=1}
describes the possible 
singular values $\gamma_1,\gamma_2$
of a $2$ by $2$ complex (or real) symmetric matrix $Z$ of the form 
$Z=\begin{smallbmatrix}
* & s_1 \\ 
s_1 & *
\end{smallbmatrix}$. 
To compare, Theorem~\ref{th:LiPoon} asserts that  the eigenvalues
$\lambda_1\geq \lambda_2$ of $Z$ satisfy 
$2s_1\leq \lambda_1-\lambda_2$ (and, generally speaking, nothing
else), which translates into the restrictions \eqref{eq:p=1} for the
singular values $\gamma_1=|\lambda_1|$ and $\gamma_2=|\lambda_2|$.  
\end{example}

\begin{example}
For $p\!=\!2$, the triples $(I,J,K)$ and the  corresponding
inequalities \eqref{eq:lr-inequality}~are: 
\begin{equation}
\label{eq:p=2}
\begin{array}{ccrl}
(\{1\},\{1\},\{1\})          & \quad &  2s_1\leq &\!\!\gamma_1+\gamma_2 \\[.05in]
(\{1\},\{2\},\{2\})          & &  2s_2\leq&\!\! \gamma_1+\gamma_4 \\[.05in]
(\{2\},\{1\},\{2\})          & &  2s_2\leq&\!\! \gamma_2+\gamma_3 \\[.05in]
(\{1,2\},\{1,2\},\{1,2\})    & & 2(s_1+s_2)\leq&\!\! \gamma_1+\gamma_2+\gamma_3+\gamma_4 
\end{array}
\end{equation}
Thus, in the special case $p=2$, $n=4$, the equivalence
\eqref{any-complex}$\Leftrightarrow$\eqref{lr-inequalities} 
means that the singular values 
$\gamma_1\geq\gamma_2\geq\gamma_3\geq\gamma_4\geq 0$
of a $4$ by $4$ complex matrix of the form 
$Z = 
\begin{smallbmatrix}
* & X \\ 
Y & *
\end{smallbmatrix}$
whose $2$ by~$2$ blocks $X$ and $Y$ have fixed singular values
$s_1\geq s_2$ satisfy 
\eqref{eq:p=2} 
and no other constraints. 
To compare, Theorem~\ref{th:LiPoon} provides the following list of
inequalities relating the eigenvalues
$\lambda_1\geq \lambda_2\geq \lambda_3\geq \lambda_4$ of 
a Hermitian 
$Z$ to the
singular values $s_1$ and~$s_2$ of $X$ and~$Y$: 
\begin{equation}
\label{eq:p=2/eigen}
\begin{array}{ccrl}
(\{1\},\{1\},\{1\})          & \quad &  2s_1\leq &\!\!\lambda_1-\lambda_4 \\[.05in]
(\{1\},\{2\},\{2\})          & &  2s_2\leq&\!\! \lambda_1-\lambda_3 \\[.05in]
(\{2\},\{1\},\{2\})          & &  2s_2\leq&\!\! \lambda_2-\lambda_4 \\[.05in]
(\{1,2\},\{1,2\},\{1,2\})    & & 2(s_1+s_2)\leq&\!\! \lambda_1+\lambda_2-\lambda_3-\lambda_4 
\end{array}
\end{equation}
To illustrate Remark~\ref{rem:cones}: for fixed singular values
$\gamma_1,\dots,\gamma_4$, there are 16 choices of signs for the
eigenvalues~$\pm\gamma_1,\dots,\pm\gamma_4$. 
For each of these choices, sorting the eigenvalues in decreasing order
produces a list $(\lambda_1,\dots,\lambda_4)$. 
The resulting 
inequalities are always (weakly) stronger
than~\eqref{eq:p=2}. (The latter corresponds to
$(\lambda_1,\lambda_2,\lambda_3,\lambda_4)
=(\gamma_1,\gamma_3,-\gamma_4,-\gamma_2)$.)
For example, if the eigenvalues are
$-\gamma_1,\gamma_2,\gamma_3,-\gamma_4$, then
$(\lambda_1,\lambda_2,\lambda_3,\lambda_4)
=(\gamma_2,\gamma_3,-\gamma_4,-\gamma_1)$, 
and \eqref{eq:p=2/eigen} yields four inequalities
\[
2s_1\leq \gamma_1+\gamma_2\,,\quad
2s_2\leq \gamma_2+\gamma_4\,,\quad
2s_2\leq \gamma_1+\gamma_3\,,\quad
2(s_1+s_2)\leq\gamma_1+\gamma_2+\gamma_3+\gamma_4\,, 
\]
which collectively imply~\eqref{eq:p=2}. 
\end{example}

\pagebreak[2]

\subsection{Singular value inequalities for arbitrary $X$ and~$Y$} 

It is natural to ask whether the restriction in 
Theorem~\ref{th:main} that $X$ and $Y$ have the same singular values 
can be removed (with some other collection of inequalities playing the role
of~\eqref{eq:lr-inequality}). 

\begin{problem}
\label{problem:XY}
Find necessary and sufficient conditions for the existence of matrices
$X$, $Y$, and $Z= \begin{smallbmatrix}* & X \\
        Y & *\end{smallbmatrix}$ 
with given singular values. 
Can those conditions be given by a collection of
linear inequalities? 
\end{problem}

In Proposition~\ref{pr:pxyq} below, we provide a set of necessary (but
not sufficient) conditions for this problem.
These conditions, however, will turn out to be sufficient in the
special case considered in Theorem~\ref{th:main},
and will play a role in the proof of the latter. 

Recall that the triples of sequences of nonnegative real numbers 
$(a_1\geq \cdots\geq a_n)$,
$(b_1\geq \cdots\geq b_n)$,
$(c_1\geq \cdots\geq c_n)$ 
that can occur as singular values of complex $n$ by $n$ matrices $A$,
$B$, and~$C=A+B$ are also characterized 
by a list of linear inequalities.
More specifically, note that the $2n$ by $2n$ matrices appearing in the identity
\[
\begin{bmatrix}
0   & A \\
A^* & 0
\end{bmatrix}
+
\begin{bmatrix}
0   & B \\
B^* & 0
\end{bmatrix}
=
\begin{bmatrix}
0   & C \\
C^* & 0
\end{bmatrix}
\]
have eigenvalues $a_1\geq\cdots\geq a_n\geq-a_n\geq\cdots\geq-a_1$, 
$b_1\geq\cdots\geq b_n\geq-b_n\geq\cdots\geq-b_1$, 
$c_1\geq\cdots\geq c_n\geq-c_n\geq\cdots\geq-c_1$, 
respectively,
and therefore these three sequences must satisfy the Horn inequalities
(see Proposition~\ref{pr:horn}).
Conversely, those inequalities are sufficient for the existence
of $A$, $B$, and $C=A+B$ with the desired singular values; see \cite{Kl2} for
the proof, and \cite[Section~5]{Fu1} for a detailed exposition. 

\begin{proposition}
\label{pr:pxyq}
Let $2p\leq n$, and let $Z = \begin{smallbmatrix}P & X \\
        Y & Q\end{smallbmatrix}$ 
be a complex matrix 
such that 
\begin{itemize}
\item[$\bullet$]
  $X$ is $p$ by $n-p$, with singular values $s_1\geq\cdots\geq s_p$; 
\item[$\bullet$]
$Y$ is $n-p$ by $p$, with singular values $t_1\geq \cdots\geq t_p$;
\item[$\bullet$]
$Z$ is $n$ by $n$, with singular values $\gamma_1\geq\cdots\geq\gamma_n$.  
\end{itemize}
Let $\sigma_1\geq\cdots\geq\sigma_{2p}$ be the decreasing
rearrangement of the numbers $s_1,\dots,s_p,t_1,\dots,t_p$. \linebreak[2]
Then, for all $1\leq m<2n$ and all triples $(E,F,G)\in\LR_m^{2n}$,  
\begin{equation} 
\label{st-general}
2\Bigl(
\sum_{\substack{g \in G\\ g\leq 2p}} \sigma_g 
-\sum_{\substack{g \in G'\\ g\leq 2p}} \sigma_g 
\Bigr)
\leq
 \sum_{\substack{e \in E\\e\leq n}} \gamma_e 
-\sum_{\substack{e \in E'\\e\leq n}} \gamma_e 
+\sum_{\substack{f \in F\\f\leq n}} \gamma_f 
-\sum_{\substack{f \in F'\\f\leq n}} \gamma_f\,, 
\end{equation}
where we denote
\[
G'=\{g\in\{1,\dots,2n\} \mid 2n+1-g\in G\},
\]
and similarly for $E$ and~$F$. 

In particular, for any triple of the form $(F,F,G)\in\LR_m^{2n}$, with
$m<2n$, 
\begin{equation} 
\label{st-general-FFG}
\sum_{\substack{g \in G\\ g\leq 2p}} \sigma_g 
-\sum_{\substack{g \in G'\\ g\leq 2p}} \sigma_g 
\leq
 \sum_{\substack{f \in F\\f\leq n}} \gamma_f 
-\sum_{\substack{f \in F'\\f\leq n}} \gamma_f\,.  
\end{equation}
\end{proposition}

\begin{proof}
In the identity
\begin{equation}
\label{eq:PXYQ+PXYQ} 
\begin{bmatrix}P & X \\ Y &
        Q\end{bmatrix}
+
\begin{bmatrix}-P & X \\ Y &
        -Q\end{bmatrix}
=
2\begin{bmatrix}0 & X \\ Y &
        0\end{bmatrix}\,,
\end{equation} 
both matrices on the left-hand side  have singular values
$\gamma_1 \geq \ldots \geq \gamma_n\,$, 
as one sees by applying \eqref{eq:block-unitary} with
        $U_1 = V_1 = \sqrt{-1}I_p$ and
$U_2 = V_2 = -\sqrt{-1}I_{n-p}\,$. 
The right-hand side 
has ordered singular values 
$2\sigma_1, \dots, 2\sigma_{2p}, 0, \ldots, 0$. 
Applying (the easy part of)  \cite[Theorem~15]{Fu1}
to the matrices in \eqref{eq:PXYQ+PXYQ}, we obtain~\eqref{st-general}.
\end{proof}



The converse of Proposition~\ref{pr:pxyq} is false.
That is, the inequalities \eqref{st-general} are necessary but
not sufficient for the existence of a matrix $Z$ with described 
properties. 
See Example~\ref{ex:p1n2}. 

\pagebreak[2]

\begin{example}
\label{ex:p1n2}
Let $p\!=\!1$ and $n\!=\!2$. 
Thus, we fix $s_1\!\geq \!0$, $t_1\!\geq\! 0$, and
$\gamma_1\!\geq\!\gamma_2\!\geq\!0$, 
and look at the matrices 
$Z\!=\!\begin{smallbmatrix}p & x \\
        y & q\end{smallbmatrix}$
with $|x|\!=\!s_1$, $|y|\!=\!t_1$, and 
with singular values $\gamma_1$ and~$\gamma_2$. \linebreak[2]
By the definition of singular values,
such a matrix exists if and only if there are complex numbers $p$ and
$q$ such that
\begin{align*}
\gamma_1^2+\gamma_2^2 &= |p|^2+|q|^2+|x|^2+|y|^2\,,\\
\gamma_1\gamma_2 &= |pq-xy|\,,
\end{align*}
for some $x$ and $y$ with $|x|=s_1$ and $|y|=t_1$. 
Equivalently,
\begin{eqnarray}
&&\exists\, p,q\in \mathbb{C}\ 
  \left\{ \begin{aligned}
           \gamma_1^2+\gamma_2^2 &= |p|^2+|q|^2+s_1^2+t_1^2\\
           \gamma_1\gamma_2 &= |pq-s_1t_1|
          \end{aligned}
  \right.\notag\\
&\Longleftrightarrow\ &
\exists\, p,q\in \mathbb{C}\ 
  \left\{ \begin{aligned}
           |p|^2+|q|^2&=\gamma_1^2+\gamma_2^2 -s_1^2-t_1^2\\
           |\gamma_1\gamma_2-s_1t_1|&\leq |pq|\leq |\gamma_1\gamma_2+s_1t_1|
          \end{aligned}
  \right.\notag\\
&\Longleftrightarrow\ &
2|\gamma_1\gamma_2-s_1t_1|\leq \gamma_1^2+\gamma_2^2 -s_1^2-t_1^2\notag\\
&\Longleftrightarrow\ &
  \left\{ \begin{aligned}
           (s_1-t_1)^2&\leq (\gamma_1-\gamma_2)^2\\
           (s_1+t_1)^2&\leq (\gamma_1+\gamma_2)^2
          \end{aligned}
  \right.\notag\\
&\Longleftrightarrow\ &
  \left\{ \begin{aligned}
           s_1+t_1&\leq \gamma_1 + \gamma_2  \\
           s_1-t_1&\leq \gamma_1 - \gamma_2 \\ 
           -s_1+t_1&\leq \gamma_1 - \gamma_2 . 
          \end{aligned}
  \right.
\label{eq:p1n2-complete}
\end{eqnarray}
Thus, in this special case, the necessary and sufficient conditions 
for the existence of such matrix~$Z$ 
are given by the linear inequalities \eqref{eq:p1n2-complete}. 
(This is precisely \cite[Lemma~1]{T2}.) 
On the other hand, 
the only essential inequalities among
\eqref{st-general} are the ones 
corresponding to $E=F=G=\{1\}$ and $E=F=G=\{1,2\}$;
they are, respectively, 
$\sigma_1 \leq\gamma_1$ and 
$\sigma_1+\sigma_2 \leq \gamma_1+\gamma_2\,$.
Equivalently, 
\begin{equation}
\label{eq:p1n2-FFG}
\begin{aligned}
\max(s_1,t_1)&\leq \gamma_1\,,\\
s_1+t_1&\leq\gamma_1+\gamma_2\,.
\end{aligned}
\end{equation}
Since the inequalities \eqref{eq:p1n2-FFG} do
not imply~\eqref{eq:p1n2-complete},  the converse 
of Proposition~\ref{pr:pxyq} fails. 
\end{example}

\pagebreak[2]

\subsection{Inequalities for Littlewood-Richardson coefficients}  

Our proof of Theorem~\ref{th:main} (see
Section~\ref{sec:proof-outline}) 
is based on Theorem~\ref{th:LiPoon}, 
Proposition~\ref{pr:pxyq}, 
and the following lemma, proved combinatorially in
Section~\ref{sec:littlewood-richardson} using a result of
Carr\'e and Leclerc~\cite{CL}.

\begin{lemma}
\label{lem:IJK-FFG}
If $(I,J,K)\in\LR_r^p$, then $(F,F,G)\in\LR_{2r}^{2p}$, 
where 
\begin{align}
\label{eq:FG-F}
F &= \{2i_1\!-\! 1, \ldots , 2i_r\! -\! 1\} \cup
\{2j_1, \ldots , 2j_r\}, \\ 
\label{eq:FG-G}
G &= \{2k_1\!-\!1, \ldots , 2k_r\!-\!1\} \cup
\{2k_1, \ldots, 2k_r\}.
\end{align}
\end{lemma}

If $\lambda$ and $\mu$ are the partitions associated to $I$ and $J$ by the
correspondence~\eqref{eq:set-to-partition}, let $\tau(\lambda,\mu)$ be
the partition corresponding to the set 
$F$ defined by~\eqref{eq:FG-F}.
The lemma says that if $\lr{\lambda,\mu}{\nu}$ is positive, then
$\lr{\tau(\lambda,\mu),\tau(\lambda,\mu)}{\tau(\nu,\nu)}$ is
also positive.  In Section~\ref{sec:littlewood-richardson} we prove
the following stronger assertion
(see Proposition~\ref{pr:lr-domination}): 

\begin{proposition}
\label{pr:lr-ineq-intro}
\emph{If $\lambda$, $\mu$, and $\nu$ satisfy $|\lambda| + |\mu| = |\nu|$, then} 
$
\lr{\lambda,\mu}{\nu} \leq
\lr{\tau(\lambda,\mu),\tau(\lambda,\mu)}{\tau(\nu,\nu)}. 
$
\end{proposition}

The correspondence 
$(\lambda,\mu) \mapsto \tau(\lambda,\mu)$ is a disguised
version of (the inverse of) the \emph{$2$-quotient} map,
well known in algebraic combinatorics of Young diagrams.
See more about this in Section~\ref{sec:littlewood-richardson}. 
A stronger but unproved inequality between Littlewood-Richardson 
coefficients (see Conjecture~\ref{conj}) is 
discussed in Section~\ref{sec:grassmann}.

\subsection{Proof of Theorem~\ref{th:main} modulo
  Theorem~\ref{th:LiPoon} and 
Proposition~\ref{pr:lr-ineq-intro}
}
\label{sec:proof-outline}

The implications 
\eqref{alternating-eigenvalues}$
\Rightarrow$\eqref{real-symmetric}$\Rightarrow$\eqref{any-complex} 
in Theorem~\ref{th:main} are trivial. 
The equivalence
\eqref{alternating-eigenvalues}$\Leftrightarrow$\eqref{lr-inequalities}
was proved in Section~\ref{sec:compare-to-LP} using Theorem~\ref{th:LiPoon}. 
It remains to prove that 
\eqref{any-complex}$\Rightarrow$\eqref{lr-inequalities}.

Let $(I,J,K)\in\LR_r^p$, and let $F$ and $G$ be given
by~\eqref{eq:FG-F}--\eqref{eq:FG-G}. 
By Lemma~\ref{lem:IJK-FFG}
(which follows from Proposition~\ref{pr:lr-ineq-intro}), 
$(F,F,G)\in\LR_{2r}^{2p}\subset \LR_{2r}^{n}$.  
We next apply Proposition~\ref{pr:pxyq}, with $t_k=s_k$ and
$\sigma_{2k-1}=\sigma_{2k}=s_k$ for
$k=1,\dots,p$. 
Observing that the negative sums on both sides of
\eqref{st-general-FFG} disappear for $F,G\subset\{1,\dots,n\}$, 
we obtain:
\[
\sum_{k \in K} s_k + \sum_{k \in K} s_k
=
\sum_{g \in G} \sigma_g 
\leq
 \sum_{f \in F} \gamma_f 
=\sum_{i \in I} \gamma_{2i-1} + \sum_{j \in J} \gamma_{2j}\,,
\]
as desired.
\qed

\begin{remark}
\label{rem:pool-together}
The argument in the proof above can be adapted to show
that condition \eqref{any-complex} of Theorem~\ref{th:main}
can be replaced by
\begin{itemize}
\item[{\rm (\ref{any-complex}${}'$)}]
There exists an $n$ by $n$ complex matrix 
$
\begin{smallbmatrix}
* & X \\
Y & * 
\end{smallbmatrix}
$
with singular values $\gamma_1,\dots,\gamma_n$ 
such that the matrix
$ \begin{smallbmatrix} 0
      & X \\ Y & 0\end{smallbmatrix}$ has singular
     values $s_1,s_1,s_2,s_2,\dots ,s_p,s_p$. 
\end{itemize}
(Here, as before, $X$ and $Y$ are $p$ by $n-p$ and $n-p$ by $p$,
respectively.) 
\end{remark}

\subsection{Outline of the rest of the paper}
In Section~\ref{sec:horn-combined}, we use Lemma~\ref{lem:IJK-FFG} to
characterize the possible eigenvalues of a matrix obtained as a sum of
two Hermitian matrices with a given combined list of eigenvalues. 
In turn, this result leads to (apparently) new statements
concerning Littlewood-Richardson coefficients. 

Sections~\ref{sec:proof-of-th:LiPoon}
and~\ref{sec:littlewood-richardson} 
contain the proofs of Theorem~\ref{th:LiPoon} and 
Proposition~\ref{pr:lr-ineq-intro}, respectively. 

Although not required for our proof, we include in
Section~\ref{sec:grassmann} a 
geometric argument that deduces Proposition~\ref{pr:lr-ineq-intro} from a
stronger inequality (see Conjecture~\ref{conj})
for the Littlewood-Richardson coefficients. 
For each ordered pair $(\lambda, \mu)$ of partitions, a
simple rule produces another pair
$(\lambda^*, \mu^*)$, with
$|\lambda^*| + |\mu^*| = |\lambda| + |\mu|$.  We conjecture
that $\lr{\lambda^* \, \mu^*}{\nu} \geq
\lr{\lambda \, \mu}{\nu}$ for all partitions~$\nu$.
While some cases of this conjecture can be deduced from known matrix
identities, the general case seems to require new ideas.

\pagebreak[3]

\section{Horn's problem for a combined list of eigenvalues}
\label{sec:horn-combined}

Throughout this section, ``Hermitian''
can be replaced by ``real symmetric.''

For any $p$-element lists 
$a = (a_1 \geq \cdots \geq a_p)$ and $b = (b_1 \geq \cdots \geq b_p)$
of real numbers, let $\Horn(a;b)$ denote
the set of $p$ by $p$ matrices $C$ that can be expressed as $C=A+B$,
where $A$ and~$B$ are Hermitian matrices with eigenvalues $a$
and~$b$, respectively. 
By the celebrated conjecture of A.~Horn's (proved by A.~A.~Klyachko
and A.~Knutson-T.~Tao, 
see~\cite[Section~1]{Fu1}), the set $\Horn(a;b)$ 
is described as follows.

\begin{proposition}
\label{pr:horn}
$\Horn(a;b)$ consists of the matrices $C$
whose eigenvalues $c_1\ge\cdots\ge c_p$ satisfy the trace condition
$
\textstyle \sum_i c_i = \sum_i a_i + \sum_i b_i \,$, 
together with the Horn inequalities 
\begin{equation}
\label{eq:horn-inequalities}
\textstyle \sum_{k \in K} c_k   \, \leq \,  
\sum_{i \in I} a_i  +  \sum_{j \in J} b_j \,, 
\end{equation}
for all $r<p$ and all triples $(I,J,K)$ in $\LR_r^p$.
\end{proposition}
  
Now suppose that rather than fixing the lists $a$ and~$b$,
we only fix their union
$\gamma=(\gamma_1\ge\cdots\ge\gamma_{2p})=a\sqcup b$ 
(taken with multiplicities). 
Which matrices $C$ can be written as a sum of two matrices
whose joint list of eigenvalues is~$\gamma$? 
According to Proposition~\ref{pr:abc'} below, the answer is given by the set
\begin{equation}
\label{eq:interlacing-cone}
\Horn(\gamma_1, \gamma_3, \dots, \gamma_{2p-1};\gamma_2, \gamma_4,
\dots, \gamma_{2p}). 
\end{equation}
In other words, any other splitting of $\gamma$ 
into two $p$-element sublists $a$ and~$b$ produces a set $\Horn(a;b)$
that is contained in \eqref{eq:interlacing-cone}.

\begin{proposition} 
\label{pr:abc'} 
Let $A$ and $B$ be $p$ by $p$ Hermitian 
matrices. 
Let $\gamma_1 \ge \cdots \ge \gamma_{2p}$ be the eigenvalues of $A$
and $B$ arranged in descending order. 
Then there exist Hermitian 
matrices $\widetilde A$ and $\widetilde B$ with eigenvalues
$\gamma_1, \gamma_3, \dots, \gamma_{2p-1}$ 
and $\gamma_2, \gamma_4, \dots, \gamma_{2p}$, respectively, such that
$\widetilde A + \widetilde B = A+B$.
\end{proposition}

Although Proposition~\ref{pr:abc'} can be deduced by comparing the equivalence 
\eqref{real-symmetric}$\Leftrightarrow$\eqref{alternating-eigenvalues}
of Theorem~\ref{th:main}
with Theorem~\ref{th:LiPoon}, we provide a shortcut proof below,
based directly on Lemma~\ref{lem:IJK-FFG}. 

\begin{proof}
Let $r \leq p$, and suppose that $(I,J,K)\in\LR_r^p$.
Let $F$ and $G$ be given by~\eqref{eq:FG-F}. 
By Lemma~\ref{lem:IJK-FFG}, $(F,F,G)\in\LR_{2r}^{2p}$.
Applying the corresponding Horn inequality to the identity
\begin{equation*}
\label{eq:AB+BA=CC} 
\begin{bmatrix}A & 0 \\ 0 & B\end{bmatrix}
+
\begin{bmatrix}B & 0 \\ 0 & A\end{bmatrix}
=
\begin{bmatrix}A+B & 0 \\ 0 &A+B\end{bmatrix}\,,
\end{equation*} 
we obtain
\[
2\Bigl(\sum_{i\in I}\gamma_{2i-1}+\sum_{j\in J} \gamma_{2j}\Bigr)
\ge 2\sum_{k\in K}c_k\,,
\] 
where $c_1,\dots,c_p$ are the eigenvalues of $C=A+B$. 
Since $\sum \gamma_{2i-1}+\sum\gamma_{2j}=\sum c_i$ as well,
the claim follows by Proposition~\ref{pr:horn}. 
\end{proof}

\begin{example}
Let $p\!=\!2$. 
We are looking at matrices $C$ that can be expressed as $A+B$,
where $A$ and $B$ have the joint list of eigenvalues 
$\gamma=(\gamma_1\ge\gamma_2\ge\gamma_3\ge\gamma_4)$. 
The eigen\-values $c_1\geq c_2$ of $C$ must satisfy
$c_1+c_2=\sum\gamma_i$, along with the
inequalities~\eqref{eq:horn-inequalities}, which depending on the
splitting of $\gamma$ into $a=(a_1\geq a_2)$ and $b=(b_1\geq b_2)$,
will take the following form: 
\begin{equation}
\label{eq:p=2/split}
\begin{array}{ccccccc}
I,J,K        &&\begin{array}{l}
              a=(\gamma_1,\gamma_3) \\ 
              b=(\gamma_2,\gamma_4) 
              \end{array}
                                    && \begin{array}{l}
                                      a=(\gamma_1,\gamma_4) \\ 
                                      b=(\gamma_2,\gamma_3) 
                                      \end{array}
                                                         &&\begin{array}{l}
                                                          a=(\gamma_1,\gamma_2) \\ 
                                                          b=(\gamma_3,\gamma_4) 
                                                          \end{array}
\\                     
\hline\\[-.17in]
\{1\},\{1\},\{1\}  &&  c_1\leq\gamma_1+\gamma_2 &&  c_1\leq\gamma_1+\gamma_2 && c_1\leq\gamma_1+\gamma_3 \\[.05in]
\{1\},\{2\},\{2\}  &&  c_2\leq\gamma_1+\gamma_4 &&  c_2\leq\gamma_1+\gamma_3 && c_2\leq\gamma_1+\gamma_4 \\[.05in]
\{2\},\{1\},\{2\}  &&  c_2\leq\gamma_2+\gamma_3 &&  c_2\leq\gamma_2+\gamma_4 && c_2\leq\gamma_2+\gamma_3 
\end{array}
\end{equation}
Replacing $c_2$ by $\gamma_1+\cdots+\gamma_4-c_1$, we obtain the
following conditions for $c_1$, for each of the three possible
splittings shown in \eqref{eq:p=2/split}:
\begin{equation*}
\begin{array}{ccccc}
              \begin{array}{l}
              a=(\gamma_1,\gamma_3) \\ 
              b=(\gamma_2,\gamma_4) 
              \end{array}
                                    && \begin{array}{l}
                                      a=(\gamma_1,\gamma_4) \\ 
                                      b=(\gamma_2,\gamma_3) 
                                      \end{array}
                                                         &&\begin{array}{l}
                                                          a=(\gamma_1,\gamma_2) \\ 
                                                          b=(\gamma_3,\gamma_4) 
                                                          \end{array}
\\                     
\hline\\[-.17in]
c_1\leq\gamma_1+\gamma_2 &&  c_1\leq\gamma_1+\gamma_2 && c_1\leq\gamma_1+\gamma_3 \\[.05in]
c_1\geq\max(\gamma_2+\gamma_3,\gamma_1+\gamma_4) 
                    &&  c_1\geq\gamma_1+\gamma_3 && c_1\geq\max(\gamma_2+\gamma_3,\gamma_1+\gamma_4) 
\end{array}
\end{equation*}
It is easy to see that the conditions in the first column are the
least restrictive among the three sets. 
To give a concrete example, take $\gamma_1=4$, $\gamma_2=3$,
$\gamma_3=2$, $\gamma_4=1$. Then
\begin{align*}
\Horn(4,2;3,1) &= \{C\,|\, c_1\in[5,7],\, c_1+c_2=10\},\\
\Horn(4,1;3,2) &= \{C\,|\, c_1\in[6,7],\, c_1+c_2=10\},\\
\Horn(4,3;2,1) &= \{C\,|\, c_1\in[5,6],\, c_1+c_2=10\}. 
\end{align*}
\end{example}

Proposition~\ref{pr:abc'} has the following direct implication. 

\begin{corollary}
\label{cor:sv-ev}
Condition~\eqref{three-hermitian} of
Corollary~\ref{cor:three-hermitian} is equivalent to
\begin{romitem}
\label{three-hermitian-combined}
There exist (positive semidefinite)
Hermitian $p$ by $p$ matrices $A$, $B$, and~$C$ such that: 
$\begin{smallbmatrix}
A & 0 \\
0 & B 
\end{smallbmatrix}$
has eigenvalues $\gamma_1,\dots,\gamma_{2p}$; 
$C$~has eigenvalues $s_1,\dots,s_p$; and
$2C\leq A+B$. 
\end{romitem}
\end{corollary}

\begin{corollary}
\label{cor:interlacing-partitions}
For a pair $(\lambda,\mu)$ of  partitions, 
let $\gamma_1 \!\ge\! \cdots\! \ge\!\gamma_{2p}$
be the decreasing rearrangement of  the $\lambda_i$ and $\mu_j$'s.
Define two partitions
\begin{equation}
\label{eq:split-even-odd}
\widetilde\lambda=(\gamma_1, \gamma_3, \dots, \gamma_{2p-1}),\quad
\widetilde\mu=(\gamma_2, \gamma_4, \dots, \gamma_{2p}).
\end{equation}
Then for every
partition $\nu$ such that $\lr{\lambda\,\mu}{\nu}>0$,  we
have $\lr{\widetilde\lambda\,\widetilde\mu}{\nu}>0$.
\end{corollary}

\begin{proof}
It is known (see~\cite{Fu1}) that 
$(I,J,K)\in\LR_r^p$ if and only if
$\lambda(I), \lambda(J), \lambda(K)$ are the eigenvalues
of some Hermitian 
matrices $X$, $Y$, and $Z\!=\!X\!+\!Y$.
The claim then follows by Proposition~\ref{pr:abc'}. 
\end{proof}

Recall that a Littlewood-Richardson coefficient $\lr{\lambda\,\mu}{\nu}$
is the coefficient of the Schur function~$s_\nu$
in the Schur function expansion of the product~$s_\lambda s_\mu$.
(For alternative representation-theoretic and intersection-theoretic
interpretations, see, e.g., \cite{Fu0,EC2}.) 
In view of this, the assertion of
Corollary~\ref{cor:interlacing-partitions} can be restated as follows.

\begin{corollary}
\label{cor:schur-domination}
Among all ways to distribute the parts $\gamma_1,\gamma_2,\dots$ 
of a given partition $\gamma=(\gamma_1\geq\gamma_2\geq\cdots)$ 
between two partitions $\lambda$ and~$\mu$, 
there is one, namely, \eqref{eq:split-even-odd}, 
for which the set of Schur functions contributing to the expansion
of~$s_\lambda s_\mu$ is largest by containment. 
\end{corollary}

\begin{conjecture}
\label{conj:schur-domination}
In Corollary~\ref{cor:interlacing-partitions},
$\lr{\widetilde\lambda\,\widetilde\mu}{\nu}\geq \lr{\lambda\,\mu}{\nu}\,$. 
Thus, in Corollary~\ref{cor:schur-domination}, 
any expression of the form $s_{\widetilde\lambda}s_{\widetilde\mu}-s_\lambda
s_\mu$ is a nonnegative linear combination of Schur functions. 
\end{conjecture}

\begin{example}
Let $\gamma=(3,2,1)$.
Then 
\begin{align}
\label{eq:31*2}
s_{\widetilde\lambda}s_{\widetilde\mu}=
s_{31} s_2 &= s_{51}+s_{42}+s_{33}+s_{411}+s_{321}\,,\\
\label{eq:32*1}
s_{32} s_1 &= s_{42}+s_{33}+s_{321}\,,\\
\label{eq:3*21}
s_3 s_{21} &= s_{51}+s_{42}+s_{411}+s_{321}\,.
\end{align}
We see that the right-hand side of \eqref{eq:31*2} dominates each of
the right-hand sides of \eqref{eq:32*1}--\eqref{eq:3*21},
in agreement with Corollary~\ref{cor:schur-domination}
and Conjecture~\ref{conj:schur-domination}. 
\end{example}

Proposition~\ref{pr:abc'} can be generalized to sums of several matrices,
as follows. 

\begin{proposition} 
\label{abc''}
Let $\gamma_1 \ge \cdots \ge \gamma_{mp}$ be the combined list of 
eigenvalues (with multi\-plicities) 
of $p$ by $p$ Hermitian matrices $A_1,\dots,A_m$. 
Then there are Hermitian matrices  $\widetilde A_1,\dots,\widetilde A_m$
such that $\sum\widetilde{A}_i =\sum A_i$, 
and each $\widetilde A_i$ has eigenvalues
$\gamma_i, \gamma_{i+m}, \dots, \gamma_{i+m(p-1)}$. 
\end{proposition}

\begin{proof}
This result is proved by combining Proposition~\ref{pr:abc'}
with an elementary combinatorial argument. 
Let us color the indexing set $\{1,\dots,mp\}$ 
according to which eigenvalue comes from which matrix. 
More precisely, we use the colors $1,\dots,m$,
each of them exactly $p$ times,
so that the following condition is satisfied:
\begin{itemize}
\item[($*$)]
for each color~$c$, 
the numbers $\gamma_j$ whose index $j$ has color~$c$
are precisely the eigenvalues of~$A_c\,$.  
\end{itemize}

Consider the following ``repainting'' operation:
pick two colors $c$ and~$c'$,
identify the $2p$ indices colored in these colors,
and repaint these indices (if needed) so that their colors interlace as
the indices increase: $c,c',c,c',\dots$. 
(Alternatively, we may repaint them $c',c,c',c,\dots$.)
The remaining $mp-2p$ indices keep their colors. 
By Proposition~\ref{pr:abc'}, this operation 
can always be combined with an appropriate change of matrices
$A_c$ and $A_{c'}$ so that condition ($*$) remains
fulfilled and the sum $\sum{A}_i$ is unchanged. 

To complete the proof, it suffices to show that any coloring of the
set $\{1,\dots,mp\}$ into colors $1,\dots,m$
(using each color $p$ times) can be transformed  
by a sequence of repainting operations 
into the canonical coloring
\[
1,\dots,m,1,\dots,m,\dots\dots,1,\dots,m,
\]
where each index $j$ has color $c$ with $j\equiv c\bmod m$. 
Suppose we have a non-canonical coloring, and let $k$ be the smallest 
index whose color differs from the canonical one. 
Say, $k$ has color~$c$, whereas in the canonical coloring, it has
color~$c'$.
By applying a repainting operation to the colors $c$ and~$c'$
we can change the color of~$k$, thus expanding the initial segment
colored in a canonical way.
Iterating this procedure, we will arrive at the canonical coloring. 
\end{proof}

Using Proposition~\ref{abc''}, one can extend
Corollary~\ref{cor:interlacing-partitions} 
to $m$-tuples of partitions.

\pagebreak[2]

\section{Proof of Theorem \ref{th:LiPoon}}
\label{sec:proof-of-th:LiPoon} 

We begin with a proposition that refines what it 
means for the Horn inequalities to hold.  It is 
essentially equivalent to the main result of~\cite{Fu2}, 
improved by an idea from~\cite{LP}.  We then use this 
proposition to give a quick proof of Theorem \ref{th:LiPoon}.

Let $a = (a_1 \geq \cdots \geq a_n)$, $b = (b_1 \geq \cdots \geq b_n)$, 
and $c = (c_1 \geq \cdots \geq c_n)$ be sequences of weakly 
decreasing real numbers of length $n$.  For every $1 \leq r \leq n$, 
we have a collection of Horn inequalities
\begin{equation}
\label{*rn}
\sum_{k \in K} c_k   \, \leq \,  
\sum_{i \in I} a_i  +  \sum_{j \in J} b_j , 
\tag{$*_r^n$}
\end{equation}
one for each $(I,J,K)$ in $\LR_r^n$.  For $r = n$, there is just one 
inequality
\begin{equation}
\label{*nn}
\sum_{k=1}^n c_k  \, \leq \,  \sum_{i=1}^n a_i  +  \sum_{j=1}^n b_j .
\tag{$*_n^n$}
\end{equation}
 
\begin{proposition}
\label{pr:HE} 
Assume that the sequences $a$, $b$, and $c$ consist of nonnegative 
real numbers.  The following conditions are equivalent:
\begin{enumerate}
\item[(i)]  The inequalities \textnormal{($*_r^n$)} are satisfied for all $r \leq n$. 
\item[(ii)]  There are Hermitian $n$ by $n$ matrices
  $A$, $B$, and $C$  
  with eigenvalues  $a_1, \ldots, a_n$, $b_1, \ldots, b_n$, and $c_1,
  \ldots, c_n$, respectively, such that $C \leq A + B$, i.e., $A + B -
  C$ is positive semidefinite.
\item[(iii)]  For some integer $s \geq 1$, 
  there are: 
  \begin{itemize} 
  \item
  real numbers $t(\ell) \in [0,1]$, for $1\leq \ell\leq s$, 
  \item 
  a decomposition $n = \sum_{\ell = 1}^s n(\ell)$ of $n$
  into a sum of $s$ positive integers, and  
  \item 
  a decomposition of each of 
  $a$, $b$, and $c$ into a union of $s$ subsequences, 
  denoted $a(\ell)$, $b(\ell)$, and $c(\ell)$, respectively,
  for $1\leq \ell\leq s$, 
  each of length $n(\ell)$,  
  \end{itemize}
such that the triples 
$$
(t(\ell)\cdot a(\ell), t(\ell)\cdot b(\ell), c(\ell))
$$
satisfy all inequalities \textnormal{($*_{r(\ell)}^{n(\ell)}$)} for $r(\ell) \leq n(\ell)$, with 
strict inequalities for $r(\ell) < n(\ell)$, and equality for $r(\ell) = n(\ell)$. 
\item[(iv)]  
For some integer $s \geq 1$, 
there are: 
  \begin{itemize} 
  \item
real numbers $t(\ell) \in [0,1]$, for $1 \leq \ell \leq s$,   
  \item
a decomposition $n = \sum_{\ell = 1}^s n(\ell)$ of $n$
into a sum of $s$ positive integers, and 
  \item
Hermitian $n(\ell)$ by $n(\ell)$ matrices 
$A(\ell)$, $B(\ell)$, and $C(\ell)$, for $1 \leq \ell \leq s$, with 
$$
C(\ell) = t(\ell)(A(\ell) + B(\ell)),
$$
such that $A(\ell)$ and $B(\ell)$ preserve no proper subspace of
$\mathbb{C}^{n(\ell)}$, and the eigenvalues of $\oplus
A(\ell)$, $\oplus B(\ell)$, and 
$\oplus C(\ell)$ are $a_1, \ldots, a_n$, $b_1, \ldots, b_n$, and $c_1,
\ldots, c_n$, respectively.
\end{itemize}
\end{enumerate}
\end{proposition}
 
\begin{proof}
The equivalence of (i) and (ii) is the main result of~\cite{Fu2} (for three matrices).  

(iii) $\Rightarrow$ (iv): By the Klyachko theorem \cite{Kl}\cite[Section 1]{Fu1}, 
there are 
Hermitian matrices $A(\ell)$, $B(\ell)$, and $C(\ell)$, with eigenvalues 
$a(\ell)$, $b(\ell)$, and $c(\ell)$, such that 
$C(\ell) = t(\ell)(A(\ell) + B(\ell))$.  If there is any proper subspace $L$ preserved by 
$A(\ell)$ and $B(\ell)$, and therefore by~$C(\ell)$, then $L^{\perp}$ is also
preserved; after changing bases one can further decompose the matrices
$A(\ell)$, $B(\ell)$, and $C(\ell)$. 

(iv) $\Rightarrow$ (ii): Take $A = \oplus A(\ell)$, $B = \oplus B(\ell)$, 
$C = \oplus C(\ell)$.  Since $A(\ell) + B(\ell)$ is positive semidefinite, $C(\ell) = t(\ell)(A(\ell) + B(\ell)) \leq A(\ell) + B(\ell)$ for each 
$\ell$, so $C \leq A + B$.

(i) $\Rightarrow$ (iii):  As in ~\cite{LP}, let $t$ be the 
smallest real number such that the triple $(t\cdot a, t\cdot b, c)$ satisfies 
all inequalities ($*_r^n$) for all $r \leq n$, but such that one 
(or more) 
of these inequalities holds with equality; let $(I,J,K) \in \LR_r^n$ be a 
triple for which equality holds.  Decompose $a$, $b$, and $c$ 
respectively into subsequences given by
\[
\begin{array}{rrr}
a' = (a_i)_{i \in I}, &
b' = (b_j)_{j \in J}, &
c' = (c_k)_{k \in K}, 
\\[.05in]
a'' = (a_i)_{i \notin I}, &
b'' = (b_j)_{j \notin J}, &
c'' = (c_k)_{k \notin K}.
\end{array}
\]
By \cite{Fu2} (see the discussion after the statement of Theorem 2), there are 
$r$ by $r$ Hermitian matrices $A'$, $B'$, and $C'$, with eigenvalues $a'$, $b'$, 
and $c'$, with $C' = t(A'+B')$, and there are $n-r$ by $n-r$ matrices 
$A''$, $B''$, and $C''$, with eigenvalues $a''$, $b''$, 
and $c''$, with $C'' \leq t(A''+B'')$. 
We may assume that $A'$ and $B'$ have no common invariant subspace, or they
could be further decomposed.   
Take $t(1) = t$, $n(1) = r$, $A(1) = A'$, $B(1) = B'$, 
and $C(1) = C'$.  The triple $(ta', tb',c')$ satisfies the 
inequalities (\ref{*rn}) strictly for $r(1) < n(1)$, since any
equality would lead to a proper invariant subspace by~\cite{Fu1},
Proposition~6. 
The inductive hypothesis applies to the triple 
$(t\cdot a'', t\cdot b'', c'')$, and this produces the other terms 
in the required decomposition. 
\end{proof}

\begin{remark} 
Even if all the inequalities (\ref{*rn}) for $r < n$ are strict for a
particular triple $(a,b,c)$, 
this does \emph{not} imply that the triple is indecomposable
(in the sense of~(iii)). 
For example, the triple $((2,1,0),(2,1,0),(3,2,1))$ satisfies all 
inequalities ($*_r^3$) strictly for $r < 3$, with equality for $r = 3$; 
but it decomposes into the three triples $((2),(1),(3))$,
$((1),(0),(1))$, and $((0),(2),(2))$.  
Note also that decompositions need not be unique, as this triple also
decomposes into $((1),(2),(3))$, $((0),(1),(1))$, and $((2),(0),(2))$.
\end{remark}

\begin{remark} 
Proposition~\ref{pr:HE} extends as usual (cf.\ \cite{Fu1, Fu2}) 
to the case where $a$ and $b$ are replaced by any number 
$m \geq 2$ of sequences, and with $m$ matrices in place of $A$ and~$B$.
Furthermore, the Hermitian 
matrices in (ii) and (iv) can be taken to be real symmetric.  (They 
may have no real invariant subspaces, even if they have complex 
invariant subspaces, but (iv) is true with either interpretation.)   
\end{remark}

We are now ready to prove Theorem~\ref{th:LiPoon}.

Let us first check that the equivalences
(ii)$\Leftrightarrow$(iii)$\Leftrightarrow$(iv) 
in Theorem~\ref{th:LiPoon} follow from~\cite{Fu2}. 
The equivalence of (ii) and (iii) is the result of~\cite{Fu2}, 
while (iii)$\Rightarrow$(iv) is obvious. 

(iv) $\Rightarrow$ (ii):  Let the eigenvalues of $A$ and $B$ be
$a_1,\ldots,a_p$ and $b_1,\ldots,b_p$.  The inequality $2C \leq A - B$
gives the inequalities 
\[
2 \textstyle\sum_{k \in K} s_k \leq \textstyle\sum_{i \in I} a_i -
\textstyle\sum_{j \in J} b_{p+1-j}
\]
for all $(I,J,K) \in \LR_r^p$, $r \leq p$.
Since $a_1\geq \dots \geq a_p$ and $ b_1 \geq
\dots \geq b_p$ are  subsequences of $\lambda_1 \geq \dots \lambda_2 \geq
\dots\geq
\lambda_n$ , it follows that $\lambda_i \geq a_i$ and
$\lambda_{n+1-i} \leq b_{p+1-i}$ for $1 \leq i \leq p$.
Hence 
\[
\textstyle\sum_{i \in I} a_i - \textstyle\sum_{j \in J} b_{p+1-j} \leq
  \textstyle\sum_{i \in I} \lambda_i - \textstyle\sum_{j \in J}
  \lambda_{n+1-j}\,,
\]
and (ii) follows. 

To prove the equivalence of (i) and (ii), we use
Proposition~\ref{pr:HE}. 
Note that both conditions
(and in fact each of (i)--(iv)) is
unchanged if every $\lambda_i$ is replaced by $\lambda_i + c$, for any
real number~$c$.  
This is obvious in~(ii), and follows by adding a
scalar matrix $cI_n$ to the matrix in~(i).
Hence we may assume that $\lambda_p \geq 0\geq \lambda_{n+1-p}$.  
Then all three sequences 
\[
(\lambda_1 \geq \dots \geq \lambda_p), \quad  
(-\lambda_n \geq\dots \geq -\lambda_{n+1-p}),  \quad  
(2s_1 \geq \ldots \geq 2s_p) 
\]
consist of nonnegative numbers, so Proposition~\ref{pr:HE} applies to them. 

(ii) $\Rightarrow$ (i):  We use a decomposition as in 
Proposition~\ref{pr:HE}(iv), but with $n$ replaced by~$p$.   \linebreak[2]
This produces a decomposition $p = \sum p(\ell)$, numbers $t(\ell) \in [0,1]$,
and $p(\ell)$ by $p(\ell)$ Hermitian matrices $A(\ell)$, $B(\ell)$,
and $C(\ell)$ such that $2C(\ell)=t(\ell)(A(\ell) - B(\ell))$, 
the eigenvalues of $\oplus A(\ell)$ are $\lambda_1, \dots, \lambda_p$,
the eigenvalues of $\oplus B(\ell)$ are $\lambda_{n+1-p},\dots,\lambda_n$, 
and the eigenvalues of $\oplus C(\ell)$ are $s_1, \dots, s_p$. 
Following~\cite{LP}, for each $\ell$, choose~$\theta(\ell)$ so that
$\sin(2\theta(\ell)) \!=\! t(\ell)$.  \linebreak[2]
Writing matrices in block form,
define $p(\ell)$ by $p(\ell)$ Hermitian (or real symmetric) matrices
$P(\ell)$ and $Q(\ell)$ by the identities
\[
\begin{bmatrix}P(\ell) & C(\ell) \\ C(\ell)^* & Q(\ell)
  \end{bmatrix} \,=\, 
\begin{bmatrix}\cos(\theta(\ell)) & -\sin(\theta(\ell)) \\
    \sin(\theta(\ell)) & \cos(\theta(\ell)) \end{bmatrix}
\cdot 
\begin{bmatrix}A(\ell) & 0 \\
0 & B(\ell) \end{bmatrix} \cdot
\begin{bmatrix}\cos(\theta(\ell)) & \sin(\theta(\ell)) 
    \\ -\sin(\theta(\ell)) & \cos(\theta(\ell))
  \end{bmatrix}.
\]
The direct sum of the matrices 
$\begin{smallbmatrix}P(\ell) & C(\ell) \\ C(\ell)^* & Q(\ell) \end{smallbmatrix}$, 
together with the diagonal matrix of size
$n-2p$ with diagonal entries $\lambda_{p+1},
\dots,\lambda_{n-p}$ has eigenvalues $\lambda_1,
\dots,\lambda_n$.  The direct sum $X$ of the matrices~$C(\ell)$,
filled in with $0$'s on the right, is the upper $p$ by $n-p$ block $X$ of
this matrix, and its singular values are $s_1, \ldots , s_p$. 

(i) $\Rightarrow$ (ii):
The argument is similar to the proof of Proposition~\ref{pr:pxyq}. 
Let $Z = \begin{smallbmatrix}P & X \\
    X^* & Q\end{smallbmatrix}$ be a matrix as in~(i).
The matrix $\widetilde{Z} = \begin{smallbmatrix}-P & X \\ X^* &
    -Q\end{smallbmatrix}$
has eigenvalues $-\lambda_n,\ldots,-\lambda_1$, 
as follows from the identity~(\ref{eq:block-unitary}),
with $U_1 = \sqrt{-1}I_p$, 
$V_1 = U_1^* = -\sqrt{-1}I_p$, $U_2 = -\sqrt{-1}I_{n-p}\,$, and 
$V_2 = U_2^* = \sqrt{-1}I_{n-p}\,$.
%
%
Furthermore, 
$Z + \widetilde{Z} =
\begin{smallbmatrix}0 & 2X \\ 2X^* & 0\end{smallbmatrix}
$
has eigenvalues 
\[
2s_1 \geq \ldots 2s_p \geq 0 \geq \ldots \geq
  0 \geq -2s_p \geq \ldots \geq -2s_1\,.
\] 
By the Horn inequalities for sums of Hermitian 
matrices (see~\cite{Fu1},\S1), we have, for all
$(I,J,K)$ in $\LR_r^p \subset \LR_r^n$, $r \leq p$, 
\[
\textstyle\sum_{k \in K} 2s_k   \leq   \sum_{i \in I} \lambda_i  +
\sum_{j \in J} (-\lambda_{n+1-j}), 
\]
which is the assertion in (ii).

That $X$ can be specified in advance follows as in
Remark~\ref{rem:specify-in-advance}. 
That the Hermitian mat\-rices can be taken to be real symmetric
follows from the analogous results in~\cite{Fu1, Fu2}. 
\qed

\pagebreak[2]

\section{$2$-quotients and Littlewood-Richardson coefficients}
\label{sec:littlewood-richardson} 

\subsection{The $2$-quotient correspondence}

The material reviewed in this section goes back to
T.~Nakayama~\cite{N} 
(in a somewhat different language). 
For detailed exposition and further references, see, e.g., 
\cite{EC2}, Exercise~7.59 and its solution, or~\cite{FS}. 

\begin{definition}
For two sets 
$I=\{i_1 < i_2 < \ldots < i_r \}$ and 
$J=\{j_1 < j_2 < \ldots < j_r \}$ of positive integers, 
define (cf.\ 
equation~\eqref{eq:FG-F}) 
\[
\tau(I,J) = \{2i_1\!-\! 1, \ldots , 2i_r\! -\! 1\} \cup
\{2j_1, \ldots , 2j_r\}\,. 
\]
It is easy to check that the corresponding map on partitions is well
defined. 
To be more precise, let $\lambda=(\lambda_1,\lambda_2,\dots)$ 
and $\mu=(\mu_1,\mu_2,\dots)$ be two integer partitions,
let $\ell(\lambda)$ (resp.,~$\ell(\mu)$) be the number of nonzero parts 
in~$\lambda$ (resp., in~$\mu$), and  
let $r\geq\max(\ell(\lambda),\ell(\mu))$. 
Then there are uniquely defined $r$-element sets of positive integers  
$I$ and $J$ that correspond to $\lambda$ and~$\mu$, respectively,
under the map~\eqref{eq:set-to-partition}. 
Furthermore, the partition $\tau(\lambda,\mu)$ that corresponds to
$\tau(I,J)$ under~\eqref{eq:set-to-partition} only depends on the partitions
$\lambda$ and $\mu$ and not on the sets $I$ and $J$ 
(that is, not on the choice of~$r$). 
If one traces the Young 
diagram of a partition by a sequence of horizontal and 
vertical steps, moving from Southwest to Northeast in a rectangle 
containing the diagrams of $\lambda$ and~$\mu$, the diagram of 
$\tau(\lambda,\mu)$ is traced, in a rectangle twice as wide in 
both directions, by 
alternating steps from $\lambda$ and $\mu$, starting with the 
first step of~$\lambda$.
 
It is easy to check that $|\tau(\lambda,\mu)| = 2(|\lambda| + |\mu|)$,
where we use the notation $|\lambda| = \sum\lambda_i\,$.
\end{definition}

\begin{example}
\label{ex:21*2}
Let $\lambda=(2,1)$ and $\mu=(2)$.
Taking $r=2$ gives $I=\{2,4\}$, $J=\{1,4\}$, $\tau(I,J)=\{2,3,7,8\}$,
and finally $\tau(\lambda,\mu)=(4,4,1,1)$. 
On the other hand, $r=3$ yields 
$I=\{1,3,5\}$, $J=\{1,2,5\}$, $\tau(I,J)=\{1,2,4,5,9,10\}$,
and again $\tau(\lambda,\mu)=(4,4,1,1)$. 

We identify each partition
$\lambda=(\lambda_1,\lambda_2,\dots)$
with the Young diagram that represents it 
(i.e., the one with row lengths
$\lambda_1,\lambda_2,\dots$). 
In this language, Example~\ref{ex:21*2} becomes 
\[
\setlength{\unitlength}{1.5pt}
\begin{picture}(20,42)(0,-40)
\put(-10,-5){\makebox(0,0){$\lambda=$}}
\put(0,0){\line(1,0){20}}
\put(0,-10){\line(1,0){20}}
\put(0,-20){\line(1,0){10}}
\put(0,0){\line(0,-1){20}}
\put(10,0){\line(0,-1){20}}
\put(20,0){\line(0,-1){10}}
\put(25,-9){\makebox(0,0){,}}
\end{picture}
\qquad\qquad\qquad
\setlength{\unitlength}{1.5pt}
\begin{picture}(20,42)(0,-40)
\put(-10,-5){\makebox(0,0){$\mu=$}}
\put(0,0){\line(1,0){20}}
\put(0,-10){\line(1,0){20}}
\put(0,0){\line(0,-1){10}}
\put(10,0){\line(0,-1){10}}
\put(20,0){\line(0,-1){10}}
\put(25,-9){\makebox(0,0){,}}
\end{picture}
\qquad\qquad\qquad\qquad
\setlength{\unitlength}{1.5pt}
\begin{picture}(20,42)(0,-40)
\put(-20,-5){\makebox(0,0){$\tau(\lambda,\mu)=$}}
\put(0,0){\line(1,0){40}}
\put(0,-10){\line(1,0){40}}
\put(0,-20){\line(1,0){40}}
\put(0,-30){\line(1,0){10}}
\put(0,-40){\line(1,0){10}}
\put(0,0){\line(0,-1){40}}
\put(10,0){\line(0,-1){40}}
\put(20,0){\line(0,-1){20}}
\put(30,0){\line(0,-1){20}}
\put(40,0){\line(0,-1){20}}
\put(45,-9){\makebox(0,0){.}}
\end{picture}
\]
\end{example}

\begin{example}
\label{ex:lambda=mu}
In the special case $\lambda=\mu=(\lambda_1,\lambda_2,\dots)$, 
one easily verifies that 
\[
\tau(\lambda,\lambda)= 
(2\lambda_1,2\lambda_1,2\lambda_2,2\lambda_2,\dots)\,. 
\]
That is, 
$\tau(\lambda,\lambda)$ is obtained from $\lambda$ by a
dilation with coefficient~$2$ (in both directions).
\end{example}

A Young diagram is called \emph{domino-decomposable} 
if it can be partitioned into disjoint
$1\times 2$ rectangles (\emph{dominoes}). 
The following result is a special case of a theorem of
T.~Nakayama~\cite{N}. 

\begin{proposition}
The Young diagram $\tau(\lambda,\mu)$ is always
domino-decomposable. 
The map $(\lambda,\mu)\mapsto\tau(\lambda,\mu)$ is a
bijection between ordered pairs of partitions (or Young diagrams), 
on one hand, and domino-decomposable Young diagrams, on another. 
\end{proposition}

The pair of partitions $(\lambda,\mu)$ that corresponds to a given
domino-decomposable Young diagram~$\tau=\tau(\lambda,\mu)$ 
is traditionally called the \emph{2-quotient} of~$\tau$. 

In the notation introduced above, Proposition~\ref{pr:lr-ineq-intro} is 
equivalent to the following:

\begin{proposition} 
\label{pr:lr-domination}
For any domino-decomposable Young diagram $\tau$ with 2-quotient $(\lambda,\mu)$,
and for any Young diagram $\nu$ with $|\nu| = |\lambda| + |\mu|$,
we have
$\lr{\lambda,\mu}{\nu} \leq 
\lr{\tau(\lambda,\mu),\tau(\lambda,\mu)}{\tau(\nu,\nu)}\,$.
\end{proposition}

\subsection{The result of Carr\'e and Leclerc} 

The proof of 
Proposition~\ref{pr:lr-domination}
is based on the version of the Littlewood-Richardson Rule due
to Carr\'e and Leclerc~\cite{CL} (see Proposition~\ref{pr:leclerc-carre} below),
which expresses a Littlewood-Richardson coefficient 
$\lr{\lambda\,\mu}{\nu}$
as the number of ``domino tableaux'' satisfying certain conditions. 
We briefly review this result here, referring the reader to \cite{CL}
or \cite{vL} for fine-print technicalities.

A (semistandard) \emph{domino tableau} $T$ of shape $\tau$ 
consists of a decomposition of $\tau$ into dominoes together with
the labelling of each domino by a positive integer.
The labelling must satisfy two conditions analogous to the usual
conditions imposed on (semistandard) Young tableaux: 
the labels weakly increase in rows 
and strictly increase in columns. 

The \emph{weight} of $T$ is the sequence $\nu=(\nu_1,\nu_2,\dots)$
in which each entry $\nu_i$ is equal to the number of labels in~$T$
equal to~$i$. 
The \emph{reading word} $w(T)$ is obtained by scanning the 
labels of $T$ column by column, right to left and top down. 
To clarify, when we read a tableau by
columns, right-to-left, an entry in a horizontal domino is skipped the
first time we trace it. 
(Carr\'e and Leclerc use the French notation, 
with the tableau flipped upside down with respect to our conventions, 
and their reading order is reverse to ours.) 

A domino tableau $T$ is called a
\emph{Yamanouchi domino tableaux} (YDT) 
if it satisfies the following additional
restriction: its reading word $w(T)$ is 
a \emph{Yamanouchi word}, or a \emph{lattice permutation} 
(see \cite[page~432]{EC2}), that is, 
\begin{align}
\label{eq:yamanouchi}
&\text{every entry $i$ appears in any
initial segment of $w(T)$ at least as many times as }\\[-.03in]
\nonumber
&\text{any entry $j>i$.}
\end{align} 
Figure~\ref{fig:YDT} lists all YDT of shape $(4,4,1,1)$, 
their respective reading words, and weights. 

\begin{figure}[h]
\begin{center}
\setlength{\unitlength}{1.5pt}
\begin{picture}(20,60)(0,-60)
\put(0,0){\line(1,0){40}}
\put(10,-20){\line(1,0){30}}
\put(0,-40){\line(1,0){10}}
\put(0,0){\line(0,-1){40}}
\put(10,-20){\line(0,-1){20}}
\put(40,0){\line(0,-1){20}}

\put(20,-10){\line(1,0){20}}
\put(0,-20){\line(1,0){10}}
\put(10,0){\line(0,-1){20}}
\put(20,0){\line(0,-1){20}}

\put(5,-10){\makebox(0,0){$1$}}
\put(15,-10){\makebox(0,0){$1$}}
\put(5,-30){\makebox(0,0){$2$}}
\put(30,-5){\makebox(0,0){$1$}}
\put(30,-15){\makebox(0,0){$2$}}

\put(20,-50){\makebox(0,0){$12112$}}

\put(20,-60){\makebox(0,0){$(3,2)$}}
\end{picture}
\qquad\qquad\qquad
\begin{picture}(20,60)(0,-60)
\put(0,0){\line(1,0){40}}
\put(10,-20){\line(1,0){30}}
\put(0,-40){\line(1,0){10}}
\put(0,0){\line(0,-1){40}}
\put(10,-20){\line(0,-1){20}}
\put(40,0){\line(0,-1){20}}

\put(0,-20){\line(1,0){10}}
\put(10,0){\line(0,-1){20}}
\put(20,0){\line(0,-1){20}}
\put(30,0){\line(0,-1){20}}

\put(5,-10){\makebox(0,0){$1$}}
\put(15,-10){\makebox(0,0){$1$}}
\put(25,-10){\makebox(0,0){$1$}}
\put(5,-30){\makebox(0,0){$2$}}
\put(35,-10){\makebox(0,0){$1$}}

\put(20,-50){\makebox(0,0){$11112$}}

\put(20,-60){\makebox(0,0){$(4,1)$}}
\end{picture}
\qquad\qquad\qquad
\begin{picture}(20,60)(0,-60)
\put(0,0){\line(1,0){40}}
\put(10,-20){\line(1,0){30}}
\put(0,-40){\line(1,0){10}}
\put(0,0){\line(0,-1){40}}
\put(10,-20){\line(0,-1){20}}
\put(40,0){\line(0,-1){20}}

\put(20,-10){\line(1,0){20}}
\put(0,-20){\line(1,0){10}}
\put(10,0){\line(0,-1){20}}
\put(20,0){\line(0,-1){20}}

\put(5,-10){\makebox(0,0){$1$}}
\put(15,-10){\makebox(0,0){$1$}}
\put(5,-30){\makebox(0,0){$3$}}
\put(30,-5){\makebox(0,0){$1$}}
\put(30,-15){\makebox(0,0){$2$}}

\put(20,-50){\makebox(0,0){$12113$}}

\put(20,-60){\makebox(0,0){$(3,1,1)$}}
\end{picture}
\qquad\qquad\qquad
\begin{picture}(20,60)(0,-60)
\put(0,0){\line(1,0){40}}
\put(10,-20){\line(1,0){30}}
\put(0,-40){\line(1,0){10}}
\put(0,0){\line(0,-1){40}}
\put(10,-20){\line(0,-1){20}}
\put(40,0){\line(0,-1){20}}

\put(20,-10){\line(1,0){20}}
\put(0,-20){\line(1,0){10}}
\put(0,-10){\line(1,0){20}}
\put(20,0){\line(0,-1){20}}

\put(10,-5){\makebox(0,0){$1$}}
\put(10,-15){\makebox(0,0){$2$}}
\put(5,-30){\makebox(0,0){$3$}}
\put(30,-5){\makebox(0,0){$1$}}
\put(30,-15){\makebox(0,0){$2$}}

\put(20,-50){\makebox(0,0){$12123$}}

\put(20,-60){\makebox(0,0){$(2,2,1)$}}
\end{picture}
\qquad\qquad\qquad
\end{center}
\caption{Yamanouchi domino tableaux of shape $(4,4,1,1)$}
\label{fig:YDT}
\end{figure}


\begin{proposition}
\label{pr:leclerc-carre}
{\rm \cite[Corollary~4.4]{CL}}
A Littlewood-Richardson coefficient $\lr{\lambda\,\mu}{\nu}$
is equal to the number of Yamanouchi domino tableaux
of shape $\tau(\lambda,\mu)$ and weight~$\nu$. 
\end{proposition} 

To illustrate, let $\lambda=(2,1)$ and $\mu=(2)$
(see Example~\ref{ex:21*2}). 
Then $\tau(\lambda,\mu)=(4,4,1,1)$. 
Comparing Proposition~\ref{pr:leclerc-carre} with 
Figure~\ref{fig:YDT}, we conclude that there are 4 nonvanishing 
Littlewood-Richardson coefficients $\lr{\lambda\,\mu}{\nu}\,$,
all equal to~$1$:
\[
 \lr{21,2}{32} =\lr{21,2}{41}
=\lr{21,2}{311}=\lr{21,2}{221}=1\,.
\]
Accordingly, the Schur functions
$s_\lambda$ and $s_\mu$ multiply as follows: 
\[
s_{21} \, s_2 = s_{32}+s_{41}+s_{311}+s_{221}\,. 
\]

We note that a direct link between Proposition~\ref{pr:leclerc-carre} and the
traditional versions of the Littlewood-Richardson Rule was established
by M.~A.~A.~van Leeuwen~\cite{vL}. 

\subsection{Proof of Proposition~\ref{pr:lr-domination}}

Set $\rho=\tau(\lambda,\mu)$. 
We need to show that $\lr{\lambda\,\mu}{\nu} \leq 
\lr{\rho\,\rho}{\tau(\nu,\nu)}$. 
By Proposition~\ref{pr:leclerc-carre}, 
these Littlewood-Richardson coefficients are given
by
\begin{align*}
\lr{\lambda\,\mu}{\nu}
&=\text{number of YDT of shape $\rho$ and
  weight~$\nu$,}\\ 
\lr{\rho,\rho}{\tau(\nu,\nu)}
&=\text{number of YDT of shape $\tau(\rho,\rho)$ and
  weight~$\tau(\nu,\nu)$.}
\end{align*}
Recall from Example~\ref{ex:lambda=mu} that $\tau(\rho,\rho)$
and~$\tau(\nu,\nu)$ 
are obtained from $\rho$ and $\nu$, respectively, by a
dilation with coefficient~$2$.
To prove the inequality, we need an injection $T\mapsto T'$ 
from the first set of YDT 
to the second one. To construct such an injection, simply chop each
domino (say, with a label~$k$) in a YDT $T$ of shape $\rho$ and
  weight~$\nu$ into $4$ quarter-size dominoes;
then put the labels $2k-1$ into the top two dominoes, and $2k$
into the bottom two. 
To illustrate, the leftmost tableau $T$ in Figure~\ref{fig:YDT}
will transform as shown in Figure~\ref{fig:YDT-doubled}.

\begin{figure}[h]
\begin{center}
\setlength{\unitlength}{3pt}
\begin{picture}(40,50)(0,-50)
\multiput(10,-0.5)(1,0){31}{\circle*{0.2}}
\multiput(10,-1.5)(1,0){31}{\circle*{0.2}}
\multiput(10,-2.5)(1,0){31}{\circle*{0.2}}
\multiput(10,-3.5)(1,0){31}{\circle*{0.2}}
\multiput(10,-4.5)(1,0){31}{\circle*{0.2}}
\multiput(10,-5.5)(1,0){31}{\circle*{0.2}}
\multiput(10,-6.5)(1,0){31}{\circle*{0.2}}
\multiput(10,-7.5)(1,0){31}{\circle*{0.2}}
\multiput(10,-8.5)(1,0){31}{\circle*{0.2}}
\multiput(10,-9.5)(1,0){31}{\circle*{0.2}}
\multiput(10,-10.5)(1,0){31}{\circle*{0.2}}
\multiput(10,-11.5)(1,0){31}{\circle*{0.2}}
\multiput(10,-12.5)(1,0){31}{\circle*{0.2}}
\multiput(10,-13.5)(1,0){31}{\circle*{0.2}}
\multiput(10,-14.5)(1,0){31}{\circle*{0.2}}
\multiput(10,-15.5)(1,0){31}{\circle*{0.2}}
\multiput(10,-16.5)(1,0){31}{\circle*{0.2}}
\multiput(10,-17.5)(1,0){31}{\circle*{0.2}}
\multiput(10,-18.5)(1,0){31}{\circle*{0.2}}
\multiput(10,-19.5)(1,0){31}{\circle*{0.2}}
\multiput(20,-20.5)(1,0){21}{\circle*{0.2}}
\multiput(20,-21.5)(1,0){21}{\circle*{0.2}}
\multiput(20,-22.5)(1,0){21}{\circle*{0.2}}
\multiput(20,-23.5)(1,0){21}{\circle*{0.2}}
\multiput(20,-24.5)(1,0){21}{\circle*{0.2}}
\multiput(20,-25.5)(1,0){21}{\circle*{0.2}}
\multiput(20,-26.5)(1,0){21}{\circle*{0.2}}
\multiput(20,-27.5)(1,0){21}{\circle*{0.2}}
\multiput(20,-28.5)(1,0){21}{\circle*{0.2}}
\multiput(20,-29.5)(1,0){21}{\circle*{0.2}}
\multiput(20,-30.5)(1,0){21}{\circle*{0.2}}
\multiput(20,-31.5)(1,0){21}{\circle*{0.2}}
\multiput(20,-32.5)(1,0){21}{\circle*{0.2}}
\multiput(20,-33.5)(1,0){21}{\circle*{0.2}}
\multiput(20,-34.5)(1,0){21}{\circle*{0.2}}
\multiput(20,-35.5)(1,0){21}{\circle*{0.2}}
\multiput(20,-36.5)(1,0){21}{\circle*{0.2}}
\multiput(20,-37.5)(1,0){21}{\circle*{0.2}}
\multiput(20,-38.5)(1,0){21}{\circle*{0.2}}
\multiput(20,-39.5)(1,0){21}{\circle*{0.2}}
\multiput(20,-40.5)(1,0){21}{\circle*{0.2}}
\multiput(20,-41.5)(1,0){21}{\circle*{0.2}}
\multiput(20,-42.5)(1,0){21}{\circle*{0.2}}
\multiput(20,-43.5)(1,0){21}{\circle*{0.2}}
\multiput(20,-44.5)(1,0){21}{\circle*{0.2}}
\multiput(20,-45.5)(1,0){21}{\circle*{0.2}}
\multiput(20,-46.5)(1,0){21}{\circle*{0.2}}
\multiput(20,-47.5)(1,0){21}{\circle*{0.2}}
\multiput(20,-48.5)(1,0){21}{\circle*{0.2}}
\multiput(20,-49.5)(1,0){21}{\circle*{0.2}}

\put(15,-30){\makebox(0,0){$k\!+\!1$}}
\put(12.5,-25){\makebox(0,0){$\bullet$}}
\put(17.5,-25){\makebox(0,0){$\bullet$}}

\thicklines
\put(0,0){\line(1,0){40}}
\put(10,-20){\line(1,0){10}}
\put(10,-40){\line(1,0){10}}

\put(10,-20){\line(0,-1){20}}
\put(20,-20){\line(0,-1){20}}


\end{picture}
\qquad\qquad
\begin{picture}(40,50)(0,-50)
\multiput(10,-0.5)(1,0){31}{\circle*{0.2}}
\multiput(10,-1.5)(1,0){31}{\circle*{0.2}}
\multiput(10,-2.5)(1,0){31}{\circle*{0.2}}
\multiput(10,-3.5)(1,0){31}{\circle*{0.2}}
\multiput(10,-4.5)(1,0){31}{\circle*{0.2}}
\multiput(10,-5.5)(1,0){31}{\circle*{0.2}}
\multiput(10,-6.5)(1,0){31}{\circle*{0.2}}
\multiput(10,-7.5)(1,0){31}{\circle*{0.2}}
\multiput(10,-8.5)(1,0){31}{\circle*{0.2}}
\multiput(10,-9.5)(1,0){31}{\circle*{0.2}}
\multiput(10,-10.5)(1,0){31}{\circle*{0.2}}
\multiput(10,-11.5)(1,0){31}{\circle*{0.2}}
\multiput(10,-12.5)(1,0){31}{\circle*{0.2}}
\multiput(10,-13.5)(1,0){31}{\circle*{0.2}}
\multiput(10,-14.5)(1,0){31}{\circle*{0.2}}
\multiput(10,-15.5)(1,0){31}{\circle*{0.2}}
\multiput(10,-16.5)(1,0){31}{\circle*{0.2}}
\multiput(10,-17.5)(1,0){31}{\circle*{0.2}}
\multiput(10,-18.5)(1,0){31}{\circle*{0.2}}
\multiput(10,-19.5)(1,0){31}{\circle*{0.2}}
\multiput(30,-20.5)(1,0){11}{\circle*{0.2}}
\multiput(30,-21.5)(1,0){11}{\circle*{0.2}}
\multiput(30,-22.5)(1,0){11}{\circle*{0.2}}
\multiput(30,-23.5)(1,0){11}{\circle*{0.2}}
\multiput(30,-24.5)(1,0){11}{\circle*{0.2}}
\multiput(30,-25.5)(1,0){11}{\circle*{0.2}}
\multiput(30,-26.5)(1,0){11}{\circle*{0.2}}
\multiput(30,-27.5)(1,0){11}{\circle*{0.2}}
\multiput(30,-28.5)(1,0){11}{\circle*{0.2}}
\multiput(30,-29.5)(1,0){11}{\circle*{0.2}}
\multiput(20,-30.5)(1,0){21}{\circle*{0.2}}
\multiput(20,-31.5)(1,0){21}{\circle*{0.2}}
\multiput(20,-32.5)(1,0){21}{\circle*{0.2}}
\multiput(20,-33.5)(1,0){21}{\circle*{0.2}}
\multiput(20,-34.5)(1,0){21}{\circle*{0.2}}
\multiput(20,-35.5)(1,0){21}{\circle*{0.2}}
\multiput(20,-36.5)(1,0){21}{\circle*{0.2}}
\multiput(20,-37.5)(1,0){21}{\circle*{0.2}}
\multiput(20,-38.5)(1,0){21}{\circle*{0.2}}
\multiput(20,-39.5)(1,0){21}{\circle*{0.2}}
\multiput(20,-40.5)(1,0){21}{\circle*{0.2}}
\multiput(20,-41.5)(1,0){21}{\circle*{0.2}}
\multiput(20,-42.5)(1,0){21}{\circle*{0.2}}
\multiput(20,-43.5)(1,0){21}{\circle*{0.2}}
\multiput(20,-44.5)(1,0){21}{\circle*{0.2}}
\multiput(20,-45.5)(1,0){21}{\circle*{0.2}}
\multiput(20,-46.5)(1,0){21}{\circle*{0.2}}
\multiput(20,-47.5)(1,0){21}{\circle*{0.2}}
\multiput(20,-48.5)(1,0){21}{\circle*{0.2}}
\multiput(20,-49.5)(1,0){21}{\circle*{0.2}}

\put(20,-25){\makebox(0,0){$k\!+\!1$}}
\put(15,-22.5){\makebox(0,0){$\bullet$}}
\put(25,-22.5){\makebox(0,0){$\bullet$}}

\thicklines
\put(0,0){\line(1,0){40}}
\put(10,-20){\line(1,0){20}}
\put(10,-30){\line(1,0){20}}

\put(10,-20){\line(0,-1){10}}
\put(30,-20){\line(0,-1){10}}


\end{picture}
\end{center}
\caption{Checking condition (iii)}
\label{fig:YDT-iii}
\end{figure}
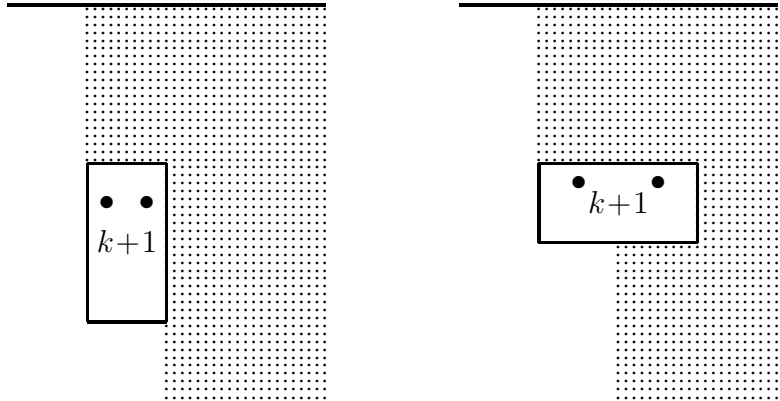

We then need to check that
\begin{itemize}
\item[(i)]  the resulting tableau $T'$ is a valid (semistandard)
  domino tableau;
\item[(ii)] $T'$ has shape $\tau(\rho,\rho)$ and
  weight~$\tau(\nu,\nu)$; 
\item[(iii)] $T'$ satisfies the Yamanouchi
  condition~\eqref{eq:yamanouchi}
for any $i$ and $j=i+1$.
\end{itemize}
Verifying the claims (i) and (ii) is straightforward. Claim (iii),
for $i$ odd, is also easy:
each entry $i+1$ is preceded by~$i$
in the reading word~$w(T')$. 
The case of $i$ even requires careful examination of a handful of cases.
For $i=2k$, we need to look at a domino labelled $k+1$ in the original
tableau $T$ (see Figure~\ref{fig:YDT-iii}) 
and check that each of the corresponding entries equal to $2k+1$ in~$T'$
(marked by a bullet~$\bullet$ in Figure~\ref{fig:YDT-iii}) 
appears in the reading word $w(T')$ at the end of an initial segment
that contains more $2k$'s than $2k+1$'s. 
This can be done by looking at all dominoes labeled $k$ or $k+1$ in the
shaded region  in Figure~\ref{fig:YDT-iii}, and invoking
condition~\eqref{eq:yamanouchi} for the tableau~$T$.
The details are left to the reader.  
\qed

\begin{figure}[h]
\begin{center}
\setlength{\unitlength}{3pt}
\begin{picture}(20,40)(0,-40)
\thicklines
\put(0,0){\line(1,0){40}}
\put(10,-20){\line(1,0){30}}
\put(0,-40){\line(1,0){10}}
\put(0,0){\line(0,-1){40}}
\put(10,-20){\line(0,-1){20}}
\put(40,0){\line(0,-1){20}}

\put(20,-10){\line(1,0){20}}
\put(0,-20){\line(1,0){10}}
\put(10,0){\line(0,-1){20}}
\put(20,0){\line(0,-1){20}}

\put(5,-10){\makebox(0,0){$1$}}
\put(15,-10){\makebox(0,0){$1$}}
\put(5,-30){\makebox(0,0){$2$}}
\put(30,-5){\makebox(0,0){$1$}}
\put(30,-15){\makebox(0,0){$2$}}

\put(50,-10){\makebox(0,0){$\longmapsto$}}
\end{picture}
\qquad\qquad\qquad\qquad\qquad 
\begin{picture}(20,40)(0,-40)
\thicklines
\put(0,0){\line(1,0){40}}
\put(10,-20){\line(1,0){30}}
\put(0,-40){\line(1,0){10}}
\put(0,0){\line(0,-1){40}}
\put(10,-20){\line(0,-1){20}}
\put(40,0){\line(0,-1){20}}

\put(20,-10){\line(1,0){20}}
\put(0,-20){\line(1,0){10}}
\put(10,0){\line(0,-1){20}}
\put(20,0){\line(0,-1){20}}

\thinlines

\put(0,-10){\line(1,0){20}}
\put(0,-30){\line(1,0){10}}
\put(20,-5){\line(1,0){20}}
\put(20,-15){\line(1,0){20}}

\put(5,0){\line(0,-1){40}}
\put(15,0){\line(0,-1){20}}
\put(30,0){\line(0,-1){20}}

\put(2.5,-5){\makebox(0,0){$1$}}
\put(7.5,-5){\makebox(0,0){$1$}}
\put(12.5,-5){\makebox(0,0){$1$}}
\put(17.5,-5){\makebox(0,0){$1$}}
\put(2.5,-15){\makebox(0,0){$2$}}
\put(7.5,-15){\makebox(0,0){$2$}}
\put(12.5,-15){\makebox(0,0){$2$}}
\put(17.5,-15){\makebox(0,0){$2$}}
\put(2.5,-25){\makebox(0,0){$3$}}
\put(7.5,-25){\makebox(0,0){$3$}}
\put(2.5,-35){\makebox(0,0){$4$}}
\put(7.5,-35){\makebox(0,0){$4$}}
\put(25,-2.5){\makebox(0,0){$1$}}
\put(35,-2.5){\makebox(0,0){$1$}}
\put(25,-7.5){\makebox(0,0){$2$}}
\put(35,-7.5){\makebox(0,0){$2$}}
\put(25,-12.5){\makebox(0,0){$3$}}
\put(35,-12.5){\makebox(0,0){$3$}}
\put(25,-17.5){\makebox(0,0){$4$}}
\put(35,-17.5){\makebox(0,0){$4$}}
\end{picture}
\qquad\qquad\qquad\qquad
\end{center}
\caption{Injection $T\mapsto T'$  in the proof of 
Proposition~\ref{pr:lr-domination}}
\label{fig:YDT-doubled}
\end{figure}
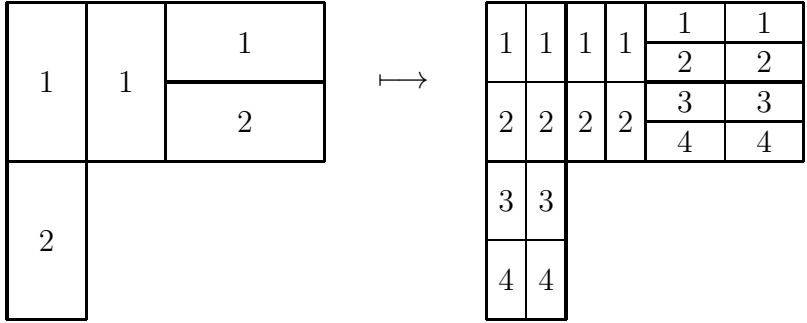

\pagebreak[2]

\section{Grassmann geometry}
\label{sec:grassmann}

In this section we sketch a geometric approach to the proof of 
Proposition~\ref{pr:lr-domination}. Carrying this out leads to
another problem about Littlewood-Richardson coefficients --- but this
remains a conjecture.  We begin by stating the conjecture, which does
not require any geometry.  

\subsection{A combinatorial conjecture}
Given an ordered pair $(\lambda, \mu)$ of partitions with the same
number of parts, define a new
ordered pair $(\lambda^*, \mu^*)$ by the following recipe:
\begin{align*}
\lambda^*_k &\eqbydef \lambda_k - k +
\#\{\ell \mid \mu_\ell - \ell \geq \lambda_k - k\}; \\
\mu^*_\ell &\eqbydef \mu_\ell - \ell + 1 +
\#\{k \mid \lambda_k - k > \mu_\ell - \ell\}.
\end{align*}
For example, if  $\lambda \!=\! (5,5,2,2)$ and $\mu \!=\!(1,1,0,0)$,
then $\lambda^* \!=\! (4,3,1,0)$ and $\mu^* \!=\! (3,2,2,1)$.

We say that a partition fits in a
$p$ by $n-p$ rectangle if it has at most $p$ positive parts, each of
size at most~$n-p$. 
Let $\lambda$ and $\mu$ be such partitions, and let
$I$ and $J$ be the $p$-element subsets of the set $\{1,\dots,n\}$ 
associated to $\lambda$ and~$\mu$, respectively, under the
correspondence \eqref{eq:set-to-partition}. Then the sets $I^*$ and
$J^*$ associated to $\lambda^*$ and~$\mu^*$ are defined by
\begin{align*}
I^* &= \bigl\{i+\#\{i'\in I\mid i'<i\}-\#\{j\in J\mid
j<i\}\bigr\}_{i\in I}\,,\\[.05in]
J^* &= \bigl\{j+\#\{j'\in J\mid j'\leq j\}-\#\{i\in I\mid
i\leq j\}\bigr\}_{j\in J}\,.
\end{align*}
Using this reformulation of the transformation $(\lambda,
\mu)\mapsto(\lambda^*, \mu^*)$, it is easy to verify that $\lambda^*$
and $\mu^*$ are partitions, that both of them fit in a $p$ by $n-p$ rectangle,
and 
$|\lambda^*| + |\mu^*| =  |\lambda| + |\mu|$.

\begin{conjecture}
\label{conj}  
For any partition $\nu$, 
we have
$\lr{\lambda \, \mu}{\nu} \leq \lr{\lambda^* \, \mu^*}{\nu}\,$. 
\end{conjecture}

Equivalently, in terms of Schur functions,
$s_{\lambda^*} s_{\mu^*} - s_{\lambda} s_{\mu}$ is
Schur positive, i.e., when this difference is expressed as a linear
combination of Schur functions $s_{\nu}$, all the coefficients are nonnegative.

Using \cite{lrcalc}, A.~Buch has verified this conjecture for all pairs
$(\lambda,\mu)$ where both $\lambda$ and $\mu$ fit in a $p$ by $q$
rectangle with $p\,q \leq 48$.

\medskip

A pair of partitions $(\lambda,\mu)$ is a fixed point of the operation
$*$ if and only if the
sequence $\mu_1, \lambda_1, \mu_2, \lambda_2, \mu_3, \ldots $ is weakly
decreasing.
In fact, for any other $(\lambda,\mu)$, if $ k = k(\lambda,\mu)$ is
maximal such that
the first $k$ terms of this sequence form a weakly decreasing sequence,
then the
corresponding sequence for $(\lambda^*,\mu^*)$ has the same first $k-1$
terms as
that for $(\lambda,\mu)$, while its $k^{\text{th}}$ term is strictly
larger; and $k(\lambda^*,\mu^*) \geq k(\lambda,\mu)$.  From this it
follows that, after applying
the $*$ operation a finite number of times, one always reaches a fixed
point --- a fact
which is also an easy consequence of the conjecture.

\subsection{Intersections of Schubert cells} 
In the rest of Section~\ref{sec:grassmann},
we present a geometric argument showing 
how Conjecture~\ref{conj} implies Proposition~\ref{pr:lr-domination}.
The general shape of the argument is as follows. 
First, we formulate a geometric conjecture
(see Conjecture~\ref{conj:transversality}) 
asserting transversality of certain intersections of Schubert cells, 
and explain why this transversality conjecture 
would imply Proposition~\ref{pr:lr-domination}. 
We then derive the transversality conjecture from
Conjecture~\ref{conj} by an analysis of
tangent spaces, using a result of P.~Belkale. 

We begin by recalling the basic facts of the Schubert calculus on
Grassmannians, while setting up the relevant notation. 
(See~\cite{Fu0} for further details.) 
Let $V$ be an $n$-dimensional vector space over an algebraically
closed field.  For any complete flag
\[
E_{\scriptscriptstyle{\bullet}}
\,=\,
(0=E_0\subset E_1\subset\cdots\subset E_n=V) 
\]
of subspaces of~$V$, and any
partition $\lambda$ whose Young diagram fits in a
$p$ by $n-p$ rectangle, there is a \emph{Schubert variety} 
$\Omega_{\lambda}(E_{\scriptscriptstyle{\bullet}})$ in the
Grassmannian $\operatorname{Gr}(p,V)$ of $p$-dimensional subspaces
of~$V$, defined by 
\[
\Omega_{\lambda}(E_{\scriptscriptstyle{\bullet}}) = 
\{ L \in \operatorname{Gr}(p,V) \mid
\dim(L \cap E_{n-p+k-\lambda_k})
\geq k \text{ for $1 \leq k \leq p$} \}.
\]
This is the closure of the corresponding \emph{Schubert cell}
$\Omega_{\lambda}^{\circ}(E_{\scriptscriptstyle{\bullet}})$, 
which consists of all subspaces $L \in \operatorname{Gr}(p,V)$ such that,
for $0\leq m\leq n$ and $0 \leq k \leq p$, we have 
$\dim(L \cap E_m)=k$ if and only if
\begin{equation}
\label{eq:schubert-cell-range}
n-p+k-\lambda_k\leq m\leq n-p+k-\lambda_{k+1}\,.
\end{equation}
If $I\subset\{1,\dots,n\}$ corresponds to $\lambda$ as
in~\eqref{eq:set-to-partition}, then
\begin{equation}
\label{eq:schubert-cell}
\Omega_{\lambda}^{\circ}(E_{\scriptscriptstyle{\bullet}})
=
\{ L \in \operatorname{Gr}(p,V) \mid
\dim(L \cap E_m)= \#\{i\in I\mid i>n-m\}
\text{ for any~$m$} \}.
\end{equation}
The Schubert cell
$\Omega_{\lambda}^{\circ}(E_{\scriptscriptstyle{\bullet}})$
is a manifold (isomorphic to an affine space) 
of codimension $|\lambda|$ in $\operatorname{Gr}(p,V)$.

Let $\lambda$, $\mu$, and $\nu$ be partitions whose Young diagrams fit
in a $p$ by $n-p$ rectangle, and assume that the Littlewood-Richardson
number
$\lr{\lambda\,\mu\,\nu}{} 
\eqbydef
\lr{\lambda \, \mu}{\nu^{\vee}}
$
is positive, where we set \linebreak[2]
$\nu^{\vee} = (n-p-\nu_p, \ldots,n-p-\nu_1)$. 
Take three complete flags $E_{\scriptscriptstyle{\bullet}}$,
$F_{\scriptscriptstyle{\bullet}}$, $G_{\scriptscriptstyle{\bullet}}$
in $V\!=\! \mathbb{C}^n$ which are in general position. 
The latter assumption implies that the corresponding Schubert cells
$\Omega_{\lambda}^{\circ}(E_{\scriptscriptstyle{\bullet}})$,
$\Omega_{\mu}^{\circ}(F_{\scriptscriptstyle{\bullet}})$,
$\Omega_{\nu}^{\circ}(G_{\scriptscriptstyle{\bullet}})$
meet transversally in $\lr{\lambda\,\mu\,\nu}{}$ points,
i.e., there are $\lr{\lambda\,\mu\,\nu}{}$
subspaces $L$ of dimension $p$ in $V$ that are in the
transversal intersection of these Schubert cells.

We next construct three complete flags
$A_{\scriptscriptstyle{\bullet}}$, $B_{\scriptscriptstyle{\bullet}}$,
and $C_{\scriptscriptstyle{\bullet}}$ in the $2n$-dimensional vector
space $V \oplus V$, by setting, for $1 \leq m \leq n$, 
\begin{equation}
\label{eq:flags-ABC}
\begin{array}{ll}
A_{2m} = F_m \oplus E_m,  &A_{2m-1} = F_{m-1}  \oplus E_m\,, \\
B_{2m} = E_m \oplus F_m,  &B_{2m-1} = E_m \oplus F_{m-1}\,, \\
C_{2m} = G_m \oplus G_m, \quad &C_{2m-1} = G_m \oplus G_{m-1}\,.
\end{array}
\end{equation}

\begin{lemma}
\label{lem:three-cells-tau}
If $L\in \Omega_{\lambda}^{\circ}(E_{\scriptscriptstyle{\bullet}})
\cap\Omega_{\mu}^{\circ}(F_{\scriptscriptstyle{\bullet}})
\cap\Omega_{\nu}^{\circ}(G_{\scriptscriptstyle{\bullet}})$, 
then
\begin{equation}
\label{eq:l-oplus-l}
L\oplus L
\in\Omega_{\tau(\lambda,\mu)}^{\circ}
(A_{\scriptscriptstyle{\bullet}})
\cap\Omega_{\tau(\lambda,\mu)}^{\circ}
(B_{\scriptscriptstyle{\bullet}})
\cap\Omega_{\tau(\nu,\nu)}^{\circ}
(C_{\scriptscriptstyle{\bullet}}).
\end{equation}
\end{lemma}

\begin{proof}
This is a straightforward verification based on the definitions 
\eqref{eq:FG-F} and~\eqref{eq:schubert-cell}; e.g.,  
\begin{align*}
 & \dim((L\oplus L)\cap A_{2m})\\
=& \dim(L\cap F_m)+\dim(L\cap E_m)\\
=& \#\{i\in I\mid i>n-m\} + \#\{j\in J\mid j>n-m\}\\
=& \#\{i\in \tau(I,J)\mid i>2n-2m\},
\end{align*}
as desired;
here $J$ and $\tau(I,J)$ denote the subsets that correspond to~$\mu$
and $\tau(\lambda,\mu)$, respectively. 
\end{proof}

\begin{conjecture}
\label{conj:transversality}
In Lemma~\ref{lem:three-cells-tau},
the intersection of Schubert cells appearing in \eqref{eq:l-oplus-l} is
transversal at each such point~$L\oplus L$.  
\end{conjecture}

\begin{remark}
The transversality would be automatic if the flags 
$A_{\scriptscriptstyle{\bullet}}$, $B_{\scriptscriptstyle{\bullet}}$,
and $C_{\scriptscriptstyle{\bullet}}$ were in general position with
respect to each other; but they are not, already for $n=1$. 
\end{remark}

Conjecture~\ref{conj:transversality}, or the weaker assertion that
these three cells  intersect properly at each such point 
(i.e., $L\oplus L$ is an
isolated point of the intersection) is enough to imply that
\[
\lr{\lambda,\,\mu,\,\nu}{} \leq 
\lr{\tau(\lambda,\mu),\,\tau(\lambda,\mu),\,\tau(\nu,\nu)}{},
\]
which is equivalent to Proposition~\ref{pr:lr-domination}. 
(Here we are using the fact that, even
if this intersection of Schubert varieties should contain connected
components of positive dimension, their contribution to the total
intersection number must be nonnegative; this follows, e.g., from the
fact that the tangent bundle of the Grassmannian is generated by its
global sections.)

\subsection{
Conjecture~\ref{conj} implies Conjecture~\ref{conj:transversality}}

Recall that the tangent space to $\operatorname{Gr}(p,V)$ at a point
$L$ is canonically identified with $\Hom(L,V/L)$. 

We use some basic facts about tangent spaces to Schubert cells, which
can be found in Belkale's preprint~\cite{B}.  If $L$ is in a Schubert cell
$\Omega_{\lambda}^{\circ}(E_{\scriptscriptstyle{\bullet}})$,
then its tangent space
$T_{[L]}(\Omega_{\lambda}^{\circ}(E_{\scriptscriptstyle{\bullet}}))$ 
is the sub\-space of $\Hom(L,V/L)$ consisting of linear 
maps $\phi$ such that 
\[
\phi(E_m \cap L) \subset (L + E_m)/L
\]
for all~$m$. 
These intersections $E_m \cap L$ and quotients $(L + E_m)/L$ give general
flags $E'_{\scriptscriptstyle{\bullet}}$ and
$E''_{\scriptscriptstyle{\bullet}}$  in $L$ and~$V/L$, respectively.  
More precisely (cf.~\eqref{eq:schubert-cell-range}),
the flag $E'_{\scriptscriptstyle{\bullet}}$ 
in $L$ is given~by
\begin{equation}
\label{eq:E-prime}
E'_k=E_{n-p+k-\lambda_k} \cap L, 
\end{equation}
for $1 \leq k \leq p$, 
and the flag $E''_{\scriptscriptstyle{\bullet}}$ 
in $V/L$ consists of the spaces $(L + E_m)/L$, 
for all $m$ not of the form $n-p+k-\lambda_k\,$.
An equivalent condition on the tangent space (which we will also use)
is that a map 
$\phi\in T_{[L]}(\Omega_{\lambda}^{\circ}(E_{\scriptscriptstyle{\bullet}}))$ 
sends $E'_k$ to $E''_{n-p-\lambda_k}$,  for $1\leq k \leq p$.
Summarizing, 
\begin{align}
\nonumber
T_{[L]}(\Omega_{\lambda}^{\circ}(E_{\scriptscriptstyle{\bullet}})) 
&= \{\phi \colon L \to V/L \mid \phi(E_m \cap L) \subset (L + E_m)/L, 
\text{ for } 1\leq m \leq n \} \\
\label{eq:tangent-space}
&=  \{\phi \colon L \to V/L \mid \phi(E'_k) \subset
E''_{n-p-\lambda_k}\,, \text{ for } 1 \leq k \leq p \}. 
\end{align}

\begin{proposition}
\label{Belkale} 
\textnormal{\cite{B}}
Assume $|\lambda| + |\mu| + |\nu| = p(n-p)$. The following are equivalent:
\begin{enumerate}
\item[(i)] $\lr{\lambda\, \mu\, \nu}{} > 0$. 
\item[(ii)] For any general complete
flags $P_{\scriptscriptstyle{\bullet}}$,
$Q_{\scriptscriptstyle{\bullet}}$, and $R_{\scriptscriptstyle{\bullet}}$
in a $p$-dimensional vector space $L$, and
$P_{\scriptscriptstyle{\bullet}}'$,
$Q_{\scriptscriptstyle{\bullet}}'$, and
$R_{\scriptscriptstyle{\bullet}}'$
in an $(n-p)$-dimensional vector space $L'$,
if $\phi \colon L \to L'$ is a linear map that sends
$P_i$ to $P_{n-p-\lambda_i}'$,
$Q_i$ to $Q_{n-p-\mu_i}'$, and
$R_i$ to $R_{n-p-\nu_i}'$, then $\phi = 0$.
\end{enumerate}
\end{proposition}

Recall that we are working under the assumptions that 
\begin{itemize}
\item
$\lr{\lambda\, \mu\, \nu}{} > 0$;
\item
$E_{\scriptscriptstyle{\bullet}}$,
$F_{\scriptscriptstyle{\bullet}}$, $G_{\scriptscriptstyle{\bullet}}$
are three complete flags in $V=\mathbb{C}^n$ in general position;
\item 
$L$ lies in the intersection of the Schubert cells
$\Omega_{\lambda}^{\circ}(E_{\scriptscriptstyle{\bullet}})$,
$\Omega_{\mu}^{\circ}(F_{\scriptscriptstyle{\bullet}})$,
$\Omega_{\nu}^{\circ}(G_{\scriptscriptstyle{\bullet}})$. 
\end{itemize}
To illustrate 
Proposition~\ref{Belkale}, let us apply it to the flags induced on $L$
and $V/L$ by $E_{\scriptscriptstyle{\bullet}}$,
$F_{\scriptscriptstyle{\bullet}}$, and $G_{\scriptscriptstyle{\bullet}}$
(cf.\ \eqref{eq:tangent-space});  
we then recover the transversality of the intersection at~$L$. 

\medskip


Recall the definition \eqref{eq:flags-ABC} of the flags 
$A_{\scriptscriptstyle{\bullet}}$, $B_{\scriptscriptstyle{\bullet}}$,
and~$C_{\scriptscriptstyle{\bullet}}$ in $V \oplus V$. 
Conjecture~\ref{conj:transversality}
says that for every $L$ as above, any linear map 
$\Phi:L \oplus L\to (V \oplus V)/(L \oplus L)$ 
that satisfies 
\begin{equation}
\label{eq:transversality-restated}
\begin{array}{l}
\Phi(A_k \cap (L \oplus L))\subset ((L \oplus L) + A_k)/(L \oplus L),\\  
\Phi(B_k \cap (L \oplus L))\subset ((L \oplus L) + B_k)/(L \oplus L),\\  
\Phi(C_k \cap (L \oplus L))\subset ((L \oplus L) + C_k)/(L \oplus L),
\end{array}
\end{equation}
for all $1\leq k\leq 2n$, must be the zero map.  Regarding 
$\Phi$ as a 2 by 2 matrix 
$\Phi=\begin{smallbmatrix} \phi_{11} & \phi_{12}\\
                           \phi_{21} & \phi_{22}\end{smallbmatrix}$
of maps from $L$ to~$V/L$, one can restate 
conditions \eqref{eq:transversality-restated} 
as saying that both $\phi_{11}$ and $\phi_{22}$ satisfy 
\begin{equation}
\label{eq:three-inclusions-diagonal}
\begin{array}{l}
\phi(E_m \cap L)\subset (L+E_m)/L,\\ 
\phi(F_m \cap L)\subset (L+F_m)/L,\\ 
\phi(G_m \cap L)\subset (L+G_m)/L, 
\end{array}
\end{equation}
while both $\phi_{12}$ and $\phi_{21}$ satisfy 
\begin{align}
\label{eq:three-inclusions-off-diagonal-1}
\phi(E_m \cap L)&\subset (L+F_m)/L,\\
\label{eq:three-inclusions-off-diagonal-2}
\phi(F_m \cap L)&\subset (L+E_{m-1})/L,\\ 
\label{eq:three-inclusions-off-diagonal-3}
\phi(G_m \cap L)&\subset (L+G_m)/L, 
\end{align}
for $1 \leq m \leq n$. 
The transversality of the three original  
Schubert cells implies that any $\phi$ satisfying
\eqref{eq:three-inclusions-diagonal} is the zero map, 
so the two diagonal maps vanish. 
It remains to show that the two off-diagonal maps vanish as well.
We will deduce it from
\eqref{eq:three-inclusions-off-diagonal-1}--\eqref{eq:three-inclusions-off-diagonal-3}
using Proposition~\ref{Belkale} and Conjecture~\ref{conj}.

Let us apply Proposition~\ref{Belkale} to the situation where
\begin{itemize}
\item
the flags
$P_{\scriptscriptstyle{\bullet}}$, $Q_{\scriptscriptstyle{\bullet}}$,
$R_{\scriptscriptstyle{\bullet}}$  are induced on $L$ 
by $E_{\scriptscriptstyle{\bullet}}$,
$F_{\scriptscriptstyle{\bullet}}$, 
$G_{\scriptscriptstyle{\bullet}}$, respectively; 
\item
$P'_{\scriptscriptstyle{\bullet}}$,  
$Q'_{\scriptscriptstyle{\bullet}}$, 
$R'_{\scriptscriptstyle{\bullet}}$ are induced on 
$L'= V/L$ by $F_{\scriptscriptstyle{\bullet}}$,
$E_{\scriptscriptstyle{\bullet}}$, 
$G_{\scriptscriptstyle{\bullet}}$, respectively (in this order!);
\item
the three partitions are $\lambda^*$, $\mu^*$, and~$\nu$.
\end{itemize}
By Conjecture~\ref{conj}, we have $\lr{\lambda^*\, \mu^*\,\nu}{} > 0$, 
so Proposition~\ref{Belkale} applies. 
We claim that unraveling condition (ii) of Proposition~\ref{Belkale}
in this particular situation, 
one obtains precisely that~\eqref{eq:three-inclusions-off-diagonal-1}--\eqref{eq:three-inclusions-off-diagonal-3}
implies $\phi=0$, as needed. 
To explain why, we will demonstrate how
\eqref{eq:three-inclusions-off-diagonal-1} 
translates into
$\phi(E'_k)\subset F''_{n-p-\lambda_k^*}$,
where we are using the same notation as in the sentence
containing~\eqref{eq:E-prime}.  
(The checks involving \eqref{eq:three-inclusions-off-diagonal-2} and
\eqref{eq:three-inclusions-off-diagonal-3} are similar.)  
It follows from the definition of a Schubert cell (see
\eqref{eq:schubert-cell-range}) that \eqref{eq:three-inclusions-off-diagonal-1} amounts to
ensuring that, for any $1\leq k\leq p$, we have
\[
\phi(E'_k)
=\phi(E_{n-p+k-\lambda_k}\cap L)
\subset(L+F_{n-p+k-\lambda_k})/L
=F''_d\,, 
\]
where $d=\dim((L+F_{n-p+k-\lambda_k})/L)$. 
It remains to calculate~$d$. 
Using the fact that $L\in
\Omega^\circ_\mu(F_{\scriptscriptstyle{\bullet}})$, 
together with~\eqref{eq:schubert-cell}, we obtain: 
\begin{align*}
d&=\dim(F_{n-p+k-\lambda_k})-\dim(F_{n-p+k-\lambda_k}\cap L)\\
&= n-p+k-\lambda_k-\#\{j\in J\mid j>\lambda_k+p-k\}\\
&= n-p+k-\lambda_k-\#\{\ell\mid \mu_\ell+p+1-\ell>\lambda_k+p-k\}\\
&= n-p-\lambda_k^*\,,
\end{align*}
as claimed. 
(In fact, this calculation is  
exactly what determines/defines $\lambda^*$ and~$\mu^*$, and  
is the only place where we need to know what they are.) 
\qed

\subsection*{Acknowledgments}

We thank Anders Buch for providing experimental evidence for 
Conjecture~\ref{conj}, and therefore for Theorem~\ref{th:main} 
before the latter was proved.  
We also thank Prakash Belkale for sending us an early version
of~\cite{B}, which motivated the reasoning in 
Section~\ref{sec:grassmann}.
Finally, we are grateful to the referee for suggesting a number of
improvements, including a simplification of the proof
of Proposition~\ref{abc''}. 

\pagebreak[2]


\end{document}